\documentclass[11pt,a4paper]{amsart}

\usepackage{amsmath, amsthm, amssymb, fullpage}

\usepackage{paralist}
\usepackage{enumerate}

\usepackage{hyperref}

\DeclareTextAccentDefault{\k}{T1}

\numberwithin{equation}{section}

\newtheorem{mainth}{Theorem}

\newtheorem{theorem}{Theorem}[section]

\newtheorem{lemma}[theorem]{Lemma}
\newtheorem{corollary}[theorem]{Corollary}
\theoremstyle{definition}
\newtheorem*{acknowledgement}{Acknowledgement}
\newtheorem{notation}[theorem]{Notation}
\theoremstyle{remark}
\newtheorem{remark}[theorem]{Remark}

\newcommand{\rd}{\mathrm{d}}
\newcommand{\cO}{\mathcal{O}}
\newcommand{\cL}{\mathcal{L}}
\newcommand{\bP}{\mathbb{P}}  
\newcommand{\bR}{\mathbb{R}}
\newcommand{\bC}{\mathbb{C}}
\newcommand{\bH}{\mathbb{H}} 
\newcommand{\bZ}{\mathbb{Z}}
\newcommand{\cm}{\operatorname{cm}}
\newcommand{\bif}{\operatorname{bif}}
\newcommand{\AP}{\mathrm{AP}}
\newcommand{\Cb}{\mathrm{Cb}}
\DeclareMathOperator{\spe}{sp}
\DeclareMathOperator{\loc}{loc} 
\newcommand{\Id}{\mathrm{Id}}
\DeclareMathOperator{\supp}{Supp}
\DeclareMathOperator{\Res}{Res}

\newcommand{\bD}{\mathbb{D}}
\newcommand{\bN}{\mathbb{N}}
\newcommand{\bQ}{\mathbb{Q}}
\DeclareMathOperator{\Per}{Per}
\DeclareMathOperator{\Int}{Int}
 
\subjclass[2020]{Primary 37F10; Secondary 31C10}
\keywords{full critically marked polynomial family,
unicritical polynomial family, multipliers, moving targets,
equidistribution, dynatomic polynomials, bifurcation measures,
thermodynamic formalism}

\begin{document}

\title{Equidistribution, multipliers, and moving targets
in polynomial families}

\author{Y\^usuke Okuyama}
\address{Division of Mathematics, Kyoto Institute of Technology, Sakyo-ku, Kyoto 606-8585 JAPAN}
\email{okuyama@kit.ac.jp}

\begin{abstract}
We study equidistribution of parameters for each of which a polynomial has a cycle
with a prescribed moving multiplier in the full critically marked or unicritical polynomial families of a given degree $d>1$.
For the full critically marked
monic polynomial family $f=(f_\lambda)$, $\lambda\in\bC^{d-1}$, of degree $d$,
letting $p_{f,n}^*(\lambda,w)$ be the multiplier polynomial for cycles
of formal exact period $n$ having multiplier $w$,
we show that for any target sequence $w_n=d^na_n$ satisfying $n^{-1}\log|a_n|\to\ell\in[-\infty,+\infty)$, the averaged current of integration
$\nu(n)^{-1}[p_{f,n}^*(\cdot,w_n)=0]$
along the zero divisor of $p_{f,n}^*(\cdot,w_n)$
tends to the $(1,1)$-current $\rd\rd^c\max\{L_d,\log d+\ell\}$
in the sense of currents on $\bC^{d-1}$ as $n\to\infty$,
where $L_d$ is the Lyapunov exponent function and $\nu(n)$ is the degree
of the $n$th formal exact dynatomic polynomial associated to $f$.
We also obtain the corresponding equidistribution theorem for
the unicritical family $f_\lambda(z)=z^d+\lambda$, $\lambda\in\bC$, of degree $d$.
The proof combines thermodynamic formalism, the Eremenko--Levin
multiplier bound, and a comparison principle for plurisubharmonic functions, in a partly quantitative manner.
A key ingredient is the vanishing of the bifurcation measure on the boundary
of every component of the interior of the connectedness locus for the full critically marked family,
which we show using Kiwi's study of critical portraits and rational laminations.
\end{abstract}

\maketitle

\section{Introduction}\label{sec:intro}

Since Brolin's seminal work \cite{Brolin} on iterated preimages of points
under polynomials, equidistribution has been one of the central viewpoints
in complex dynamics.
The study of its parameter-space analog was initiated by Levin
\cite{Levin90} for monic centered unicritical polynomial families.
Equidistribution in parameter spaces has subsequently been extensively
studied by Bassanelli--Berteloot
\cite{BassanelliBerteloot09,BassanelliBerteloot11},
Buff--Gauthier \cite{BG14}, and Gauthier \cite{Gauthier13},
and, in a quantitative way, by Gauthier--Vigny \cite{GauthierVigny17},
Okuyama \cite{Okuyama18}, and Gauthier--Okuyama--Vigny
\cite{GauthierOkuyamaVigny19}.
Our aim in this paper is to contribute to the study of equidistribution
in polynomial families of parameters for each of which a polynomial
has a cycle of period $n$ having multiplier $w_n$,
for any given sequence $(w_n)$ satisfying a logarithmic growth condition,
that is, for moving multiplier targets, as $n\to\infty$.

\subsection{Critically marked polynomials}\label{sec:crit}

Let $f(\lambda,z)=f_\lambda(z)\in\bQ[c_1,\ldots,c_{d-2},a,z]$ 
be the full critically marked monic polynomial family of degree $d$
on $\bP^1$, where $d\ge 2$ is fixed,
parametrized by
$\lambda=(c_1,\ldots,c_{d-2},a)\in\bC^{d-1}$ and determined by the condition
\begin{gather*}
 f_\lambda'(z)=d\,z\prod_{j=1}^{d-2}(z-c_j)\quad\text{and}\quad f_\lambda(0)=a^d
\end{gather*}
on the derivative and the constant term, where $\prod_{j=1}^0=1$ as usual, 
so that for every $\lambda$, $c_0:=0,c_1,\ldots,c_{d-2}$ are all the
marked $d-1$ critical points in $\bC$ of $f_\lambda$ taking into account
their multiplicities.
The Lyapunov exponent function for $f$ on the parameter space $\bC^{d-1}$ is defined as
\begin{gather*}
 L_d:\bC^{d-1}\ni\lambda\mapsto L(f_\lambda)\in[\log d,+\infty),
\end{gather*}
where $L(f_\lambda)$ is the Lyapunov exponent of $f_\lambda$ with respect to the
measure $\mu_{f_\lambda}$ of maximal entropy, and is indeed a 
proper, continuous, and plurisubharmonic function on $\bC^{d-1}$ 
(for the range and the properness, see \eqref{eq:MP} and \cite[Proposition 6.3]{DF08}), 
and the bifurcation current and measure 
for $f$ are respectively the positive closed $(1,1)$-current and the positive measure
\begin{gather*}
 T_d:=\rd\rd^c L_d\quad\text{and}\quad\mu_d:=T_d^{d-1}\quad\text{on }\bC^{d-1}.
\end{gather*}
The connectedness and escaping loci of $f$ are respectively defined as
\begin{align*}
 C_d&:=\bigl\{\lambda\in\bC^{d-1}:(f_\lambda^n(c_j))_{n\ge 0}
 \text{ is bounded for every }j\in\{0,\ldots,d-2\}\bigr\},\\
 E_d&:=\bigl\{\lambda\in\bC^{d-1}:\lim_{n\to\infty}f_\lambda^n(c_j)=\infty
 \text{ for every }j\in\{0,\ldots,d-2\}\bigr\}\subset\bC^{d-1}\setminus C_d,
\end{align*}
which are a compact connected subset in $\bC^{d-1}$
(see \cite[Remark 12.1(b), p. 320]{BrannerHubbard92})
and a domain in $\bC^{d-1}$, respectively.
The bifurcation measure $\mu_d$ is non-trivial and of finite mass,
and its support coincides with the Shilov boundary of $C_d$
(\cite[Proposition 7]{DF08}).

Let $\mu:\bN\to\{0,\pm 1\}$ denote the M\"obius function.
For each $n\in\bN$, the $n$th formal exact dynatomic polynomial associated to $f$ is 
\begin{gather*}
  \Phi_{f,n}^*(\lambda,z)
 :=\prod_{m\mid n}\bigl(f_\lambda^m(z)-z\bigr)^{\mu(n/m)}\in\bQ[c_1,\ldots,c_{d-2},a,z],
\end{gather*}
which is monic in $(\bQ[\lambda])[z]$.
For every $\lambda\in\bC^{d-1}$, a zero in $\bC$ of $\Phi_{f,n}^*(\lambda,\cdot)$ 
is called a periodic point of $f_\lambda$
having the formal exact period $n$, and the symbol * in the notation $\Phi^*_{f,n}$
means this formal exactness; for more details, 
see \cite{MortonSilverman95}. For each $n\in\bN$, setting
\begin{gather*}
  \nu(n):=\deg_z\Phi_{f,n}^*,
\end{gather*}
which is $\sim d^n$ as $n\to\infty$, 
we have $d^n=\sum_{m\mid n}\nu(m)$, and
the $n$th multiplier polynomial
$p_{f,n}^*(\lambda,w)\in\bQ[c_1,\ldots,c_{d-2},a,w]$
associated to $f$ is the unique polynomial whose leading coefficient in 
$(\bQ[\lambda])[w]$ 
is $(-1)^{\nu(n)/n}$ such that
\begin{gather*}
p_{f,n}^*(\lambda,w)^n
  =\Res_z\bigl(
    \Phi_{f,n}^*(\lambda,z),
    (f_\lambda^n)'(z)-w
  \bigr);
\end{gather*}
for more details, see \cite{MortonVivaldi95}. 
For each $n\in\bN$ and each $w\in\bC$, the zero divisor
\begin{gather*}
 \Per_f^*(n,w):=[p_{f,n}^*(\cdot,w)=0]
\end{gather*}
of $p_{f,n}^*(\cdot,w)$ is regarded as a $(1,1)$-current of integration
on the parameter space $\bC^{d-1}$.
Here and below, we use the normalization
$\rd\rd^c=(i/\pi)\partial\bar\partial$.

Our first principal result is the following.

\begin{mainth}\label{th:main}
Let $d\ge 2$, let $(a_n)_{n\in\bN}$ be a sequence in $\bC$, set $w_n:=d^na_n$
for each $n\in\bN$,
and assume that under the convention $\log 0=-\infty$, the limit
$\ell:=\lim_{n\to\infty}\log(|a_n|^{1/n})\in[-\infty,+\infty)$
exists. Then 
\begin{equation}\label{eq:main}
 \lim_{n\to\infty}\frac{1}{\nu(n)}\Per_f^*(n,w_n)
 =
 \rd\rd^c\max\{L_d,\log d+\ell\}
\end{equation}
as $(1,1)$-currents on $\bC^{d-1}$.
In particular, if $\ell\le 0$, then
\begin{gather*}
  \lim_{n\to\infty}\frac{1}{\nu(n)}\Per_f^*(n,w_n)=T_d
\end{gather*}
as $(1,1)$-currents on $\bC^{d-1}$.
\end{mainth}

The case of $w_n\equiv w$ for some $w\in\bC$ specializes to
the theorem of Gauthier \cite[Theorem 1]{Gauthier13}.
For holomorphic families of rational functions of degree $>1$,
the case $w_n\equiv w$ for $|w|<1$ was established
by Bassanelli--Berteloot
\cite{BassanelliBerteloot09,BassanelliBerteloot11};
for the full critically marked polynomial family here, they also treated $|w|=1$
\cite[Theorem 1.2]{BassanelliBerteloot11}.
The strategy in the proof of Theorem \ref{th:main} is similar to that in
\cite[Theorem 1]{Gauthier13}: we propagate
our initial analysis on $E_d$ to $\bC^{d-1}\setminus C_d$,
and finally to the whole $\bC^{d-1}$.
Our initial estimates on $E_d$ are quantitative.

Our second principal result is the following,
which plays a key role in the final propagation step above
and is derived from the finiteness result in Kiwi's study on rational laminations and
critical portraits parameterwise,
reaping it in the parameter space $\bC^{d-1}$
using Dujardin--Favre \cite{DF08}.

\begin{mainth}\label{th:boundary}\label{th:conn}
For every component $U$ of $\Int C_d$, we have $\mu_d(\partial U)=0$.
\end{mainth}

This vanishing statement in Theorem \ref{th:conn} is the same as
that in \cite[Theorem 4]{Gauthier13}\footnote{One point in \cite[\S 5.2, proof of Theorem 5.5]{Gauthier13}:
After the application of
\cite[Theorem 5.14]{McMullen94}, the nontriviality of the limit in measure
appears to require further justification, since the supports of the
invariant line fields have areas tending to zero under the contradiction
hypothesis.}.
We will mention the unicritical polynomial version of Theorem \ref{th:boundary}
in Subsection \ref{sec:unicrit}. Our argument
also yields, a bit surprisingly,
a much simpler proof of this version than the previously known proofs
which will also be mentioned
in Subsection \ref{sec:unicrit}.

\subsection{Unicritical polynomials}\label{sec:unicrit}

Let $f(\lambda,z)=f_\lambda(z):=z^d+\lambda$, $\lambda\in\bC$, be the monic and centered unicritical degree $d$ polynomial family on $\bP^1$ parametrized by $\lambda$, where $d\ge 2$ is fixed. The compact, connected, regular, and full subset 
\begin{gather*}
  M_d:=\{\lambda\in\bC: (f_\lambda^n(0))_{n\ge 0}
  \text{ is bounded}\}
\end{gather*}
in $\bC$ is still called the connectedness locus of the family $f$, and we set
$I_d:=\bC\setminus M_d$ and $B_d:=\partial M_d$. Indeed $I_d=\{\lambda\in\bC:\lim_{n\to\infty}f_\lambda^n(0)=\infty\}$, the Julia set $J_{f_\lambda}$
of $f_\lambda$ is either connected or totally disconnected, and the latter case occurs 
if and only if $\lambda\in I_d$. For each $\lambda\in\bC$, 
the escaping rate function
\begin{gather*}
  g_{f_\lambda}(z):=\lim_{m\to\infty}\frac{1}{d^m}\log^+|f_\lambda^m(z)|\quad\text{on }\bC
\end{gather*}
of $f_\lambda$ is indeed the extended Green function of $J_{f_\lambda}$ with pole $\infty$,
and gives the extended Green function $g_{I_d}$ on $I_d$ with pole $\infty$ so that
\begin{equation}\label{eq:green}
  g_{I_d}(\lambda)=g_{f_\lambda}(\lambda)=d\cdot g_{f_\lambda}(0)\quad\text{on }\bC.
\end{equation}
The bifurcation measure for the family $f$ is
\begin{gather*}
  \mu_{B_d}:=\rd\rd^c g_{I_d}.
\end{gather*}
For more details, see, e.g., the book \cite{BM01}.

For each $n\in\bN$, let $\Phi_{f,n}^*\in\bZ[\lambda,z]$ and in turn
$p_{f,n}^*\in\bZ[\lambda,w]$ be the $n$th formal exact dynatomic and multiplier
polynomials associated to $f$, respectively, which are defined
in the same manner as that in the full critically marked polynomial family case, using the M\"obius inversion
and the resultant.
For each $n\in\bN$ and each $w\in\bC$, we have
\begin{equation}\label{eq:degree}
  \deg_\lambda p_{f,n}^*(\lambda,w)
  \equiv\frac{d-1}{d}\nu(n),
\end{equation}
so that $\Per_f^*(n,w)$
is regarded as a purely atomic measure on $\bC$
of total mass $(d-1)\nu(n)/d$.

Our argument is adapted to this unicritical polynomial family $f_\lambda(z)=z^d+\lambda$
by taking a finite holomorphic self-covering $\lambda=a^d$ of the parameter space $\bC$.

\begin{mainth}\label{th:unicrit}
Let $d\ge 2$, let $(a_n)_{n\in\bN}$ be a sequence in $\bC$, set $w_n:=d^na_n$ for each $n\in\bN$,
and assume that under the convention $\log 0=-\infty$, the limit
$\ell:=\lim_{n\to\infty}\log(|a_n|^{1/n})\in[-\infty,+\infty)$
exists. Then 
\begin{equation}\label{eq:mainuni}
 \lim_{n\to\infty}\frac{d}{(d-1)\nu(n)}\Per_f^*(n,w_n)
 =
 \rd\rd^c\max\Bigl\{g_{I_d},\frac{d}{d-1}\max\{0,\ell\}\Bigr\}
\end{equation}
weakly on $\bC$.
In particular, if $\ell\le 0$, then
\begin{gather*}
  \lim_{n\to\infty}\frac{d}{(d-1)\nu(n)}\Per_f^*(n,w_n)=\mu_{B_d}
\end{gather*}
weakly on $\bC$.
\end{mainth}

The case of $w_n\equiv 0$ goes back to Levin's seminal work \cite[\S\S 2.3--2.6, pp. 3513--3514]{Levin90}.
For $d=2$, Theorem \ref{th:unicrit} specializes 
to the theorem of Buff--Gauthier \cite{BG14}\footnote{One point in \cite[\S 2]{BG14}:
In their exterior coordinate $\zeta$, where  
$\lambda=\phi_{M_d}^{-1}(\zeta^2)$ and $\phi_{M_d}:I_d\to\{|\zeta|>1\}$ is the normalized
conformal representation of $I_d$, they pointed out the commutativity
$\sigma_n(\omega\zeta)=\omega\sigma_n(\zeta)$ on $|\zeta|>1$ 
for every $\omega^n=1$ and every $n$ and used it to deduce $\lim_{n\to\infty}\sigma_n(\zeta)=2\zeta$. See also \cite{FG21}.}. 
For the unicritical polynomial family, the assertion corresponding to Theorem \ref{th:boundary} is
$\mu_{B_d}(\partial U)=0$ for every component $U$ of $\Int(M_d)$. For $d=2$, in \cite[Proof of Lemma 2]{BG14}, 
this vanishing is derived from Zakeri \cite{Zakeri06}.  
For general $d$, this vanishing follows from Graczyk--{\'S}wi\c{a}tek
\cite[Theorem 1.2 and Corollaries 1--2]{GS00} and
\cite[Theorem 1, Fact 5.2, and Proposition 8]{GS17}.

In Section \ref{sec:prelims}, we collect the potential-theoretic and dynamical preliminaries
and state the thermodynamic estimates used below.
In Section \ref{sec:complex}, we show Theorem \ref{th:main}
and explain how to adapt the argument to show Theorem \ref{th:unicrit}.
The holomorphic disks constructed in Lemma \ref{th:stretch} are also used in
Section \ref{sec:boundary}, where we show Theorem \ref{th:boundary} for the full critically marked polynomial family
using the combinatorics of external rays in the phase space $\bC$
of each polynomial.
We show the thermodynamic estimates in Appendix \ref{th:thermo}
and present a qualitative and ergodic proof of Lemma \ref{th:upper}
in Appendix \ref{sec:ergodic}.

\section{Preliminaries}
\label{sec:prelims}

\subsection{Potential-theoretic and dynamical preliminaries}
\label{sec:prelim}

Let $f$ be the full critically marked monic polynomial family of degree $d$
defined in Subsection \ref{sec:crit}. We denote by $\|\cdot\|$
the Hermitian norm on $\bC^{d-1}$.

\begin{remark}
 The family $f$ is same as that used in \cite[\S 2.1]{Gauthier13} up to complex scalings,
 and the parameter-space results in \cite{Gauthier13} are applied to $f$.
\end{remark}

For every $\lambda\in\bC^{d-1}$, we write
$J_{f_\lambda}$ for the Julia set of $f_\lambda$. The escaping rate function
\begin{gather*}
  g_{f_\lambda}(z):=\lim_{m\to\infty}\frac{1}{d^m}\log^+|f_\lambda^m(z)|\quad\text{on }\bC
\end{gather*}
of $f_\lambda$ is indeed the extended Green function of the Julia set $J_{f_\lambda}$ 
of $f_\lambda$ with pole $\infty$.
Let $\mu_{f_\lambda}$ be the equilibrium measure of $f_\lambda$ on $\bC$, 
which is indeed the unique measure of maximal entropy for $(f_\lambda,J_{f_\lambda})$ \cite{Brolin} and computed as $\mu_{f_\lambda}=\rd\rd^c g_{f_\lambda}$.
Recall that, for every $\lambda\in\bC^{d-1}$, the Lyapunov exponent of $f_\lambda$
with respect to $\mu_{f_\lambda}$ is defined and computed as
\begin{gather}
  L_d(\lambda)=L(f_\lambda):=\int_{\bC}\log|f_\lambda'|\,\mu_{f_\lambda}
  =\log d+\sum_{j=0}^{d-2}g_{f_\lambda}(c_j),\quad\text{so that}\label{eq:MP}\\
 E_d=\bigcap_{0\le j\le d-2}\{\lambda\in\bC^{d-1}:g_{f_\lambda}(c_j)>0\}
\quad\text{and}\quad
 C_d=\{\lambda\in\bC^{d-1}:L_d(\lambda)=\log d\}\notag
\end{gather}
(see \cite{Manning84,Przytycki85}). Our initial analysis in the proof of
Theorem \ref{th:main} is on the escaping locus $E_d$,
and the following, due to Eremenko--Levin, is one of the key facts in working there.

\begin{theorem}[{\cite[Theorem 1.6, (1.9)]{EL97}}]\label{th:EL}
For every $\lambda\in E_d$, every $n\in\bN$,
and every fixed point $z\in\bC$ of $f_\lambda^n$,
\begin{gather}\label{eq:EL}
 \frac{\log\bigl|(f_\lambda^n)'(z)\bigr|}{n}
 \ge (d-1)\min_{0\le j\le d-2}g_{f_\lambda}(c_j)>0.
\end{gather}
\end{theorem}
The following version of comparison principle for plurisubharmonic functions by Gauthier
is very useful. For the one-dimensional version, see Buff--Gauthier \cite[Lemma 3]{BG14}. 

\begin{lemma}[{\cite[Theorem 2]{Gauthier13}}]\label{th:compare}
Let $X$ be a complex manifold of dimension $\ge 1$ admitting a smooth
plurisubharmonic function $\phi$ such that $(\rd\rd^c\phi)^{\dim X}$ on $X$ is a
non-trivial volume form outside a proper analytic subset of $X$.
Let $\Omega\subset X$ be a domain with $C^1$ boundary, 
$K\subset\Omega$ a compact subset, and $v$ a plurisubharmonic function on
$\Omega$ such that $v$ is continuous, that $(\rd\rd^cv)^{\dim X}$ is of finite mass
and supported on $\partial K$, that $(\rd\rd^cv)^{\dim X}(\partial U)=0$ for every component $U$ of $\Int K$. Then for every plurisubharmonic function $u$ on $\Omega$
such that $u\le v$ on $\Omega$ and $u=v$ on $\Omega\setminus K$, 
we have $u=v$ on $\Omega$.
\end{lemma}

\subsection{Thermodynamic preliminaries}\label{sec:thermo}
For every $\lambda\in E_d$, the Julia set $J_{f_\lambda}$ is a Cantor
set and $(f_\lambda,J_{f_\lambda})$ is expanding locally uniformly for $\lambda\in E_d$, 
so in particular
$(f_\lambda)_{\lambda\in E_d}$ is $J$-stable (see Ma\~n\'e--Sad--Sullivan \cite{MSS}).
For each $t\in\bC$, the transfer operator with the weight $-t\log|f_\lambda'|$
on the complex Banach space $C^0(J_{f_\lambda})$ is
defined by
\begin{equation}\label{eq:ruelle}
  (\cL_{\lambda,t}h)(z)
  :=\sum_{w\in f_\lambda^{-1}(z)}\bigl|f_\lambda'(w)\bigr|^{-t}h(w),
  \quad z\in J_{f_\lambda}, \quad\text{for each }h\in C^0(J_{f_\lambda}).
\end{equation}
By a thermodynamic argument (see e.g.\ \cite{Zins00}), 
for each $(\lambda,t)\in E_d\times\bR$,
the $\bC$-linear operator $\cL_{\lambda,t}$
has a positive simple eigenvalue
$\Lambda_{\lambda,t}$,
and has the unique strictly positive eigenfunction 
$h_{\lambda,t}\in C^0(J_{f_\lambda})$
and the unique probability eigenmeasure $m_{\lambda,t}$ on $J_{f_\lambda}$
of $\cL_{\lambda,t}$ and the transpose ${}^t\cL_{\lambda,t}$, respectively, 
belonging to the (common) eigenvalue $\Lambda_{\lambda,t}$ 
normalized as $m_{\lambda,t}(h_{\lambda,t})=1$; using the strict positivity of $h_{\lambda,t}$,
such $\Lambda_{\lambda,t}$ is unique, and at $t=0$,
\begin{gather*}
 \Lambda_{\lambda,0}=d,\quad h_{\lambda,0}=1,\quad\text{and}\quad
m_{\lambda,0}=\mu_{f_\lambda}.
\end{gather*}
Fixing any compact subset $E$ in $E_d$, also by a stability argument, the partial derivative
of the (local) function $(\lambda,t)\mapsto\log\Lambda_{\lambda,t}$ in $t$ is computed as
\begin{gather*}
  \partial_t\log\Lambda_{\lambda,t}
   =-\int_{J_{f_\lambda}}\log|f_\lambda'|\,h_{\lambda,t}\,m_{\lambda,t},
\end{gather*}
which is even jointly continuous in $(\lambda,t)\in E\times[-\tau,\tau]$
for $0<\tau\ll 1$, and at $t=0$,
\begin{gather}
  \partial_t\log\Lambda_{\lambda,t}\bigm|_{t=0}
   =-L_d(\lambda)
 \Bigl(=-\log d-\sum_{j=0}^{d-2}g_{f_\lambda}(c_j)\quad\text{by }\eqref{eq:MP}\Bigr).\label{eq:deriv} 
\end{gather}
Moreover, for every $\lambda\in E$ and every $t\in[-\tau,\tau]$, 
\begin{gather}\label{eq:sum}
\int_{\bC}\bigl|(f_\lambda^n)'\bigr|^{-t}[f_\lambda^n-\Id=0]
=O(\Lambda_{\lambda,t}^n)\quad\text{as }n\to\infty,
\end{gather}
where the implicit constant in $O$ is independent
of $\lambda$ and $t$.
For more details, see Appendix \ref{th:thermo}.

\section{Proofs of Theorems \ref{th:main} and \ref{th:unicrit}}\label{sec:complex}

Fix a sequence $(w_n=d^na_n)_{n\in\bN}$ in $\bC$ as in Theorem \ref{th:main},
so that 
\begin{gather*}
 \ell=\log\lim_{n\to\infty}|a_n|^{1/n}\in[-\infty,+\infty). 
\end{gather*}
Let us show 
\begin{gather}
  \lim_{n\to\infty}\frac{1}{\nu(n)}\log|p_{f,n}^*(\cdot,w_n)|
  =\max\{L_d,\log d+\ell\}
\label{eq:full}
\end{gather}
in $L^1_{\loc}(\bC^{d-1})$,
which together with a continuity of $\rd\rd^c$ will complete the proof of Theorem \ref{th:main}. 
Let us first focus on the open superlevel set $\Omega_d(\ell)$ of $L_d$,
where for every $c\in[-\infty,+\infty)$, we set
\begin{gather}
 \Omega_d(c):=
\Bigl\{\lambda\in\bC^{d-1}:L_d(\lambda)>\log d+\max\{0,c\}\Bigr\}.\label{eq:loverl}
\end{gather}
By the case of $w_n\equiv 0$ in Theorem \ref{th:main},
which is done by Bassanelli--Berteloot, we have
\begin{gather}
  \lim_{n\to\infty}\frac{1}{\nu(n)}\log|p_{f,n}^*(\cdot,0)|
  =L_d
  \quad\text{in }L^1_{\loc}(E_d\cap\Omega_d(\ell)).\label{eq:levin}
\end{gather}
For every $n\in\bN$, the function $\log|p_{f,n}^*(\cdot,0)|$ is pluriharmonic on $E_d$ 
since $\supp\Per^*_f(n,0)\cap E_d=\emptyset$, and in turn,
for every $n\in\bN$ and every $w\in\bC$, the difference
$\log|p_{f,n}^*(\lambda,w)/p_{f,n}^*(\lambda,0)|$ is plurisubharmonic on $E_d$;
indeed,
\begin{gather} 
\frac{1}{\nu(n)}\log|p_{f,n}^*(\lambda,w)|
=\frac{1}{n}\frac{1}{\nu(n)}\int_{\bC}\log|(f_\lambda^n)'(\cdot)-w|[\Phi^*_{f,n}(\lambda,\cdot)=0],\quad\text{so that}\label{eq:part}\\
 \frac{1}{\nu(n)}\log\Biggl|\frac{p_{f,n}^*(\lambda,w)}{p_{f,n}^*(\lambda,0)}\Biggr|
=\frac{1}{n}\frac{1}{\nu(n)}\int_{\bC}\log\biggl|1-\frac{w}{(f_\lambda^n)'(\cdot)}\biggr|[\Phi^*_{f,n}(\lambda,\cdot)=0].\label{eq:diff}
\end{gather}
We show the following key lemma in a quantitative form for our future study,
using the thermodynamic input from Subsection \ref{sec:thermo}.

\begin{lemma}\label{th:upper} 
Locally uniformly on $E_d\cap\Omega_d(\ell)$,
\begin{gather*}
\frac{1}{\nu(n)}\log\biggl|\frac{p_{f,n}^*(\lambda,w_n)}{p_{f,n}^*(\lambda,0)}\biggr|\le o(1)\quad\text{as }n\to\infty.
\end{gather*}
\end{lemma}

\begin{proof}
Noting that $\log|1-z|\le \log(1+|z|)\le |z|^t/t$ for every 
$t\in(0,1]$ and recalling the joint continuity of $\partial_t\log\Lambda_{\lambda,t}$ in Subsection \ref{sec:thermo}, on every compact subset 
\begin{gather*}
 E(\epsilon):=\Bigl\{\lambda\in\bC^{d-1}:\|\lambda\|\le \frac{1}{\epsilon},\quad
 \min_{0\le j\le d-2}g_{f_\lambda}(c_j)\ge \epsilon,\quad\text{and}\quad
 L_d(\lambda)\ge \log d+\epsilon+\max\{0,\ell\}\Bigr\}
\end{gather*}
in $E_d\cap\Omega_d(\ell)$ exhausting $E_d\cap\Omega_d(\ell)$ as $0<\epsilon\searrow0$, for $0<t\ll 1$, we compute
\begin{align*}
  &\frac{1}{\nu(n)}\log\biggl|\frac{p_{f,n}^*(\lambda,w_n)}{p_{f,n}^*(\lambda,0)}\biggr|\\
&\le \frac{1}{tn}\frac{1}{\nu(n)}
    \exp\Bigl(t(n\log d+\log|a_n|)+\log\Lambda_{\lambda,t}^n+O(1)\Bigr)
& (\text{by }\eqref{eq:sum})\\
  &=\frac{1}{tn}\frac{d^n\cdot O(1)}{\nu(n)}
    \exp\biggl(n\Bigl(t(\log d+\frac{\log|a_n|}{n})
      +\log\Lambda_{\lambda,t}-\log\Lambda_{\lambda,0}\Bigr)\biggr) & (\text{by }\Lambda_{\lambda,0}=d)\\
  &\le \frac{1}{tn}\frac{d^n\cdot O(1)}{\nu(n)}
    \exp\biggl(nt\Bigl(\log d+\frac{\log|a_n|}{n}
      +\partial_s\log\Lambda_{\lambda,s}\Bigr|_{s=0}+\frac{\epsilon}{2}\Bigr)\biggr)\\
  &=\frac{1}{tn}\frac{d^n\cdot O(1)}{\nu(n)}
    \exp\biggl(nt\Bigl(\frac{\log|a_n|}{n}+\log d
      -L_d(\lambda)+\frac{\epsilon}{2}\Bigr)\biggr)& (\text{by }\eqref{eq:deriv})\\
  &\le \frac{1}{tn}\frac{d^n\cdot O(1)}{\nu(n)}
    \exp\biggl(nt\Bigl(\frac{\log|a_n|}{n}+\frac{\epsilon}{2}
-(\epsilon+\max\{0,\ell\})\Bigr)\biggr)& (\text{by the choice of }E(\epsilon))\\
&\le O(1)\cdot\frac{e^{-\frac{\epsilon}{4}tn}}{tn}\to 0\quad\text{as }n\to\infty& (\text{by }\nu(n)\sim d^n),
\end{align*}
where, for fixed $\epsilon>0$, the implicit constant in $O$ is independent of
$\lambda\in E(\epsilon)$ and $0<t\ll 1$. Now we are done by fixing any such $0<t\ll 1$.
\end{proof}

We also show the following in a quantitative form using the Eremenko--Levin lower 
bound.

\begin{lemma}\label{th:lower}
Locally uniformly on the open superlevel set 
\begin{gather*}
 \Omega'_d(\ell)=\biggl\{\lambda\in\bC^{d-1}:(d-1)\min_{0\le j\le d-2}g_{f_\lambda}(c_j)
>\log d+\max\{0,\ell\}(\ge \lim_{n\to\infty}\log|w_n|^{1/n})\biggr\}
\end{gather*}
of $(d-1)\min_{0\le j\le d-2}g_{f_\lambda}(c_j)$, so that $\Omega'_d(\ell)\subset\Omega_d(\ell)\cap E_d$ by \eqref{eq:MP}, we have
\begin{gather*}
  \frac{1}{\nu(n)}\log\Bigl|\frac{p_{f,n}^*(\lambda,w_n)}{p_{f,n}^*(\lambda,0)}\Bigr|\ge o(1)\quad\text{as }n\to\infty.
\end{gather*}
\end{lemma}

\begin{proof}
For every $\lambda\in\Omega'_d(\ell)$, if $n\gg 1$, then 
using \eqref{eq:diff} and \eqref{eq:EL}, we compute
\begin{multline*}
  \frac{1}{\nu(n)}
 \log\biggl|\frac{p_{f,n}^*(\lambda,w_n)}{p_{f,n}^*(\lambda,0)}\biggr|\\
  \ge \frac{1}{n}\frac{1}{\nu(n)}
  \int_{\bC}
  \log\biggl(1-\biggl(\frac{|w_n|^{1/n}}{\exp((d-1)\min_{0\le j\le d-2}g_{f_\lambda}(c_j))}\biggr)^n\biggr)[\Phi_{f,n}^*(\lambda,\cdot)=0]\\
\ge \frac{1}{n}
\log\biggl(1-\exp\Bigl(
-\frac{1}{2}n\Bigl((d-1)\min_{0\le j\le d-2}g_{f_\lambda}(c_j)-(\log d+\max\{0,\ell\})\Bigr)\Bigr)\biggr)\to 0
\end{multline*}
as $n\to\infty$, locally uniformly on $\Omega'_d(\ell)$ since the continuous function
$(d-1)\min_{0\le j\le d-2}g_{f_\lambda}(c_j)-(\log d+\max\{0,\ell\})$ is strictly positive there.
\end{proof}

The following is based on a standard radial quasiconformal deformation,
and the holomorphic disk $\lambda(t)$, which we construct explicitly in the proof,
also plays a key role in the proof of Theorem \ref{th:boundary}.
Set $\bH_+:=\{\Re t>0\}$.

\begin{lemma}\label{th:stretch}
For every $\lambda\in\bC^{d-1}$, 
there exists a holomorphic disk $\lambda(t)\in(\cO_{\bH_+})^{d-1}$ in $\bC^{d-1}$
passing through $\lambda$ at $t=1$ such that on $\bH_+$,
\begin{gather}\label{eq:stretch}
 g_{f_{\lambda(t)}}\bigl(c_j(\lambda(t))\bigr)
 =(\Re t)g_{f_\lambda}(c_j)\text{ for every }j\in\{0,\ldots,d-2\}\quad\text{and}\\
 L_d(\lambda(t))=(\Re t)L_d(\lambda)+(1-\Re t)\log d.\label{eq:sliding}
\end{gather}
\end{lemma}

\begin{proof}
Pick $\lambda=(c_1,\ldots,c_{d-2},a)\in\bC^{d-1}$.
The B\"ottcher coordinate for the polynomial $f_\lambda\in\bC[z]$ 
is the unique conformal isomorphism 
\begin{gather*}
 \phi_{f_\lambda}:\{g_{f_\lambda}>\max_{0\le j\le d-2}g_{f_\lambda}(c_j)\}
\to\{|\zeta|>e^{\max_{0\le j\le d-2}g_{f_\lambda}(c_j)}\}
\end{gather*}
normalized as $\lim_{z\to\infty}\phi_{f_\lambda}(z)/z=1$, and then
\begin{gather*}
 \phi_{f_\lambda}\circ f_\lambda=\phi_{f_\lambda}^d\quad\text{and}\quad
g_{f_\lambda}=\log|\phi_{f_\lambda}|,\quad\text{so that}\\
2\partial_zg_{f_\lambda}=\frac{1}{d^n}\frac{(\phi_{f_\lambda}\circ f_\lambda^n)'}{\phi_{f_\lambda}\circ f_\lambda^n}
=\frac{1}{d^n}\frac{(\phi_{f_\lambda}'\circ f_\lambda^n)\cdot (f_\lambda^n)'}{\phi_{f_\lambda}\circ f_\lambda^n}\quad\text{for }n\gg 1\text{ locally uniformly on }\{g_{f_\lambda}>0\}.
\end{gather*}
For each $t\in\bH_+$, the injective quasiconformal map
\begin{gather*}
 H_t(z):=\phi_{f_\lambda}(z)|\phi_{f_\lambda}(z)|^{t-1}=z|z|^{t-1}\bigl(1+o(1)\bigr)
 \quad\text{as }z\to\infty 
\end{gather*}
from $\{g_{f_\lambda}>\max_{0\le j\le d-2}g_{f_\lambda}(c_j)\}\cup\{\infty\}$ into $\bP^1$ has the 
Beltrami differential
\begin{gather*}
 \mu^{H_t}:=\frac{\partial_{\bar z}H_t}{\partial_zH_t}
 \frac{\rd\bar z}{\rd z}
 =\frac{(t-1)\phi_{f_\lambda}^2|\phi_{f_\lambda}|^{t-3}
 \overline{\phi_{f_\lambda}'}}
 {(t+1)|\phi_{f_\lambda}|^{t-1}\phi_{f_\lambda}'}
 \frac{\rd\bar z}{\rd z}
 =\frac{t-1}{t+1}
 \frac{(\overline{\phi_{f_\lambda}'/\phi_{f_\lambda}})}
 {\phi_{f_\lambda}'/\phi_{f_\lambda}}
 \frac{\rd\bar z}{\rd z}
 =\frac{t-1}{t+1}
 \frac{(\overline{\partial_zg_{f_\lambda}})}{\partial_zg_{f_\lambda}}
 \frac{\rd\bar z}{\rd z},
\end{gather*}
which yields the Beltrami differential 
\begin{gather*}
 \mu_{H_t}:=\begin{cases}
	     0 & \text{on }\{g_{f_\lambda}\le 0\}\\
	     (f_\lambda^n)^*\mu^{H_t},\ n\gg 1, &\text{on }\{g_{f_\lambda}>0\}  
	    \end{cases}\quad\text{on }\bP^1.
\end{gather*}
Let $h_t$ be the quasiconformal homeomorphic endomorphism of $\bP^1$ fixing $0,1,\infty$
having the Beltrami differential $\mu_{h_t}=\mu_{H_t}$, so that $h_1=\Id_{\bP^1}$ and 
that $h_t\circ f_\lambda\circ h_t^{-1}\in \cO_{\bH_+}[z]$
of degree $d$ having the $d+1$ coefficients holomorphic on $\bH_+$
(see, e.g., \cite[Theorem 4.30, the corollary to Theorem 4.37, and p. 111]{ImayoshiTaniguchi92}). Replacing $h_t$ with $b(t)\cdot h_t$ if necessary,
where $b(t)\in \cO_{\bH_+}$ is the $(d-1)$th root of 
the leading coefficient of $h_t\circ f_\lambda\circ h_t^{-1}\in \cO_{\bH_+}[z]$
such that $b(1)=1$,
and then choosing a $d$th root $a(t)\in \cO_{\bH_+}$ of $h_t(a^d)$,
which either $\equiv 0$ or vanishes nowhere on $\bH_+$, such that $a(1)=a$, 
and setting 
\begin{gather*}
 \lambda(t):=\bigl(h_t(c_1),\ldots,h_t(c_{d-2}),a(t)\bigr)\in(\cO_{\bH_+})^{d-1},
\end{gather*}
we have
\begin{gather*}
f_{\lambda(t)}'(z)=(h_t\circ f_\lambda\circ h_t^{-1})'(z)
 =d\,z\prod_{j=1}^{d-2}\bigl(z-h_t(c_j)\bigr)\quad\text{and}\quad
 f_{\lambda(t)}(0)=a(t)^d\quad\text{for every }t\in\bH_+.
\end{gather*}
Now noting that $H_t\circ h_t^{-1}$ is conformal near $\infty$
and fixes $\infty$, we have
\begin{gather*}
 \log|h_t(w)|=\log|H_t(w)|+O(1)
 =(\Re t)\log|w|+O(1)\quad\text{as }w\to\infty,\quad\text{so that}\\
 g_{f_{\lambda(t)}}(h_t(z))
 =\lim_{n\to\infty}\frac{1}{d^n}
 \log^+\bigl|f_{\lambda(t)}^n(h_t(z))\bigr|
 =\lim_{n\to\infty}\frac{1}{d^n}
 \log^+\bigl|h_t(f_\lambda^n(z))\bigr|
 =(\Re t)g_{f_\lambda}(z)\quad\text{on }\bC,
\end{gather*}
so in particular \eqref{eq:stretch} holds, and in turn
so does \eqref{eq:sliding} by \eqref{eq:MP}.
\end{proof}

Recall that $\Omega_d(c)$ is the 
open superlevel set of $L_d$ defined in \eqref{eq:loverl}.
The following is an improvement of \cite[\S 6.2]{Gauthier13}.

\begin{lemma}\label{th:extend}
For every $c\in[-\infty,+\infty)$ and every plurisubharmonic function $u$ on $\bC^{d-1}$
such that $u=L_d$ on $E_d\cap\Omega_d(c)$,
we have $u=L_d$ on $\Omega_d(c)$.
\end{lemma}

\begin{proof}
For $0\le k\le d-2$, set
\begin{gather*}
 U_k:=\bigl\{\lambda\in\Omega_d(c):
 \#\{j\in\{0,\ldots,d-2\}:g_{f_\lambda}(c_j)=0\}\le k\bigr\},
 \quad\text{so that}\quad U_{d-2}=\Omega_d(c).
\end{gather*}
Pick a plurisubharmonic function $u$ on $\bC^{d-1}$ such that 
$u=L_d$ on $U_0=E_d\cap\Omega_d(c)$, and 
let us show $u=L_d$ on $U_k$ for every $0\le k\le d-2$.
The case of $k=0$ is clear. Suppose that $1\le k\le d-2$ and
that $u=L_d$ on $U_{k-1}$. 

First, for every distinct $i,j\in\{0,\ldots,d-2\}$, using
the B\"ottcher coordinates $\phi_{f_\lambda}$ from the proof of Lemma \ref{th:stretch},
locally on
$\{\lambda\in\bC^{d-1}:g_{f_\lambda}(c_i)>0\}\cap\{\lambda\in\bC^{d-1}:g_{f_\lambda}(c_j)>0\}$, for $n\gg 1$,
\begin{gather*}
 g_{f_\lambda}\bigl(c_i(\lambda)\bigr)-g_{f_\lambda}\bigl(c_j(\lambda)\bigr)
 =\frac{1}{d^n}\log\biggl|
 \frac{\phi_{f_\lambda}\bigl(f_\lambda^n(c_i(\lambda))\bigr)}
 {\phi_{f_\lambda}\bigl(f_\lambda^n(c_j(\lambda))\bigr)}
 \biggr|
\end{gather*}
does not identically vanish by a standard quasiconformal perturbation
near (the backward orbit of) one of the two critical points, based on
\cite[\S 3, Lemma 1]{Shishikura},
so that the subset $E_{i,j}:=\{\lambda\in\bC^{d-1}:g_{f_\lambda}(c_i)=g_{f_\lambda}(c_j)>0\}$ in $\bC^{d-1}$ has real $2d-2$ Lebesgue measure $0$.

Next, pick $\lambda_0\in U_k\setminus (U_{k-1}\cup\bigcup_{i\neq j}E_{i,j})$ and set
\begin{gather*}
 I_0:=\bigl\{j\in\{0,\ldots,d-2\}:
 g_{f_{\lambda_0}}(c_j(\lambda_0))=0\bigr\},\quad \#I_0=k.
\end{gather*}
For $n\gg1$, using the vanishing of the self-intersection of each activity current
\cite[Proposition 6.9]{DF08} and \eqref{eq:MP}
(see also \cite[proof of Lemma 4.5]{Gauthier13}), pick a $k$-dimensional irreducible component $C_0$ passing 
through $\lambda_0$ of 
\begin{multline*}
\Bigl\{
\lambda\in\bC^{d-1}:
\bigl(\phi_{f_\lambda}\bigl(f_\lambda^n(c_j(\lambda))\bigr)\bigr)_{j\in\{0,\ldots,d-2\}\setminus I_0}
= 
\bigl(\phi_{f_{\lambda_0}}\bigl(f_{\lambda_0}^n(c_j(\lambda_0))\bigr)\bigr)_{j\in\{0,\ldots,d-2\}\setminus I_0}
\quad\text{and}\\
 \max_{0\le j\le d-2}g_{f_\lambda}(c_j(\lambda))<
d^n\min_{j\in\{0,\ldots,d-2\}\setminus I_0}g_{f_\lambda}(c_j(\lambda))
\Bigr\}
\end{multline*}
and then set
\begin{gather*}
 Y_0:=C_0\cap \Bigl\{\lambda\in\bC^{d-1}:
 d\max_{j\in I_0}g_{f_\lambda}(c_j(\lambda))
 <\min_{j\in\{0,\ldots,d-2\}\setminus I_0}
 g_{f_{\lambda_0}}(c_j(\lambda_0))\Bigr\}
\end{gather*}
and let $p_0:X_0\to Y_0$ be a resolution of singularities (neither $Y_0$ nor $X_0$ is necessarily connected), 
so that on $X_0$,
\begin{gather*}
 g_{f_{p_0(\cdot)}}(c_j\circ p_0)
 =g_{f_{\lambda_0}}(c_j(\lambda_0))>0
 \quad\text{for every }j\in\{0,\ldots,d-2\}\setminus I_0,\quad\text{so in particular}
\end{gather*}
\begin{multline*}
   L_d\circ p_0
 =\log d+\sum_{j\in\{0,\ldots,d-2\}\setminus I_0}g_{f_{\lambda_0}}(c_j(\lambda_0))+\sum_{j\in I_0}g_{f_{p_0(\cdot)}}(c_j\circ p_0)\\
= L_d(\lambda_0)+\sum_{j\in I_0}g_{f_{p_0(\cdot)}}(c_j\circ p_0)
\ge L_d(\lambda_0)>\log d+\max\{0,c\},
\end{multline*}
that is, $X_0\subset p_0^{-1}(\Omega_d(c))$ (cf.\, \cite[proof of Lemma 4.5]{Gauthier13}).
Hence by the induction assumption on $k-1$, we have
\begin{gather*}
 u\circ p_0=L_d\circ p_0
 \quad\text{on }\{L_d\circ p_0>L_d(\lambda_0)\}\subset p_0^{-1}(U_{k-1}).
\end{gather*}
Let us see that for $0<\epsilon\ll1$,
\begin{gather*}
 u\circ p_0\le L_d\circ p_0
 \quad\text{on }\{L_d\circ p_0<L_d(\lambda_0)+\epsilon\},
\end{gather*}
which will conclude $u\circ p_0\le L_d(\lambda_0)$ on $\{L_d\circ p_0\le L_d(\lambda_0)\}$;
indeed, for $0<\epsilon\ll1$, 
on the closure of 
$\{\lambda\in Y_0:L_d(\lambda)\le L_d(\lambda_0)+\epsilon\}$ in $\bC^{d-1}$, 
using \eqref{eq:MP}, we have
\begin{gather*}
 d\max_{j\in I_0}g_{f_\lambda}\bigl(c_j(\lambda)\bigr)
 \le d\bigl(L_d-L_d(\lambda_0)\bigr)
 \le d\epsilon
 <\min_{j\in\{0,\ldots,d-2\}\setminus I_0}g_{f_{\lambda_0}}\bigl(c_j(\lambda_0)\bigr),\quad\text{and then}
\end{gather*}
\begin{multline*}
  \max_{0\le j\le d-2}g_{f_\lambda}\bigl(c_j(\lambda)\bigr)
 =\max_{j\in\{0,\ldots,d-2\}\setminus I_0}g_{f_{\lambda}}\bigl(c_j(\lambda)\bigr)
 =\max_{j\in\{0,\ldots,d-2\}}g_{f_{\lambda_0}}\bigl(c_j(\lambda_0)\bigr)\\
 <d^n\min_{j\in\{0,\ldots,d-2\}\setminus I_0}g_{f_{\lambda_0}}\bigl(c_j(\lambda_0)\bigr)
 =d^n\min_{j\in\{0,\ldots,d-2\}\setminus I_0}g_{f_\lambda}\bigl(c_j(\lambda)\bigr)
\end{multline*}
by the choice of $n$, so that 
$\{\lambda\in Y_0:L_d(\lambda)\le L_d(\lambda_0)+\epsilon\}$ is closed in $\bC^{d-1}$,
and in turn is compact in $Y_0$ by the properness of $L_d$.
Hence using the properness of $p_0$, by the maximum principle for 
the subharmonic function $u\circ p_0$, which is $\equiv L_d\circ p_0$ 
on $\partial\{L_d\circ p_0< L_d(\lambda_0)+\epsilon\}$, 
we are done.

Once the inequality $u\circ p_0\le L_d(\lambda_0)$ on the compact subset
$\{L_d\circ p_0\le L_d(\lambda_0)\}$ in $X_0$ is at our disposal,
we follow \cite[\S 6.2]{Gauthier13}.
By
\cite[Proposition 6.9]{DF08} again and \eqref{eq:MP},
\begin{gather*}
 (\rd\rd^c(L_d\circ p_0))^k
 =k!\bigwedge_{j\in I_0}\rd\rd^c
 \bigl(g_{f_{p_0(\cdot)}}(c_j\circ p_0)\bigr)\quad\text{on }X_0.
\end{gather*}
For every $j\in I_0$, the function $g_{f_{p_0(\cdot)}}(c_j\circ p_0)$
is pluriharmonic on
$\{g_{f_{p_0(\cdot)}}(c_j\circ p_0)>0\}$ and vanishes on
$\Int\{L_d\circ p_0=L_d(\lambda_0)\}$,
so $(\rd\rd^c(L_d\circ p_0))^k$ is supported on
$\partial\{L_d\circ p_0=L_d(\lambda_0)\}$, and is of finite mass
by the continuity of $L_d\circ p_0$ and the compactness of $\{L_d\circ p_0=L_d(\lambda_0)\}$.
Moreover, by \cite[Theorem 4.2(4)]{Gauthier13},
$(\rd\rd^c(L_d\circ p_0))^k$ assigns mass $0$ to the
boundary of each component of $\Int\{L_d\circ p_0=L_d(\lambda_0)\}$,
so by Lemma \ref{th:compare} applied to $u\circ p_0$, $L_d\circ p_0$, and
$\phi=\|p_0\|^2$ on each component $X$ of a relatively compact open subset in $X_0$ 
such that this open subset contains $\{L_d\circ p_0=L_d(\lambda_0)\}$
and $X$ has $C^1$ boundary,
we in particular have $u(\lambda_0)=L_d(\lambda_0)$. 
Consequently, together with the induction assumption,
$u=L_d$ on $U_k\setminus\bigcup_{i\neq j}E_{i,j}$, and in turn
on the whole $U_k$ by plurisubharmonicity, and we are done.
\end{proof}

By the locally uniformly boundedness of the coefficients of $f_\lambda$
as holomorphic functions in $\lambda\in\bC^{d-1}$
and \eqref{eq:part}, the family
\begin{gather*}
 \biggl\{\frac{1}{\nu(n)}\log|p_{f,n}^*(\cdot,w_n)|\biggr\}_{n\in\bN}
\end{gather*}
of plurisubharmonic functions on $\bC^{d-1}$ is locally uniformly upper bounded,
and in turn using the convergence to $L_d$ in $L^1_{\loc}(\Omega'_d(\ell))$
(by \eqref{eq:levin} and Lemmas \ref{th:upper} and \ref{th:lower}) 
and H\"ormander's compactness principle for plurisubharmonic functions
(\cite[Theorem 4.1.9(a)]{Hormander83}),
pick any such plurisubharmonic function 
\begin{align*}
 u&:=\lim_{j\to\infty}\frac{1}{\nu(n_j)}\log|p_{f,n_j}^*(\cdot,w_{n_j})|
\quad\text{in }L^1_{\loc}(\bC^{d-1})
\end{align*}
on $\bC^{d-1}$
so that $u=L_d$ on $\Omega'_d(\ell)$.
By \eqref{eq:levin} and Lemma \ref{th:upper}, we have
$u\le L_d$ on Lebesgue a.e.\ $E_d\cap\Omega_d(\ell)$, and in turn, 
by the submean value inequality for the subharmonic function $u$ and
the continuity of $L_d$,
$u\le L_d$ on the whole $E_d\cap\Omega_d(\ell)$.
For each $\lambda\in E_d\cap\Omega_d(\ell)$, 
let $\lambda(t):\bH_+\to\bC^{d-1}$ be the holomorphic disk in $\bC^{d-1}$
passing through $\lambda(1)=\lambda$ obtained in Lemma \ref{th:stretch}, so that
by \eqref{eq:stretch} and \eqref{eq:sliding},
the restriction of $\lambda(t)$ to
$\{t\in\bC:(\Re t)(L_d(\lambda)-\log d)>\max\{0,\ell\}\}\subset\bH_+$
is a holomorphic disk in $E_d\cap\Omega_d(\ell)$ and
even a holomorphic disk in $\Omega'_d(\ell)$ if $\Re t\gg 1$. Then the subharmonic function
\begin{gather*}
 u\bigl(\lambda(t)\bigr)-L_d\bigl(\lambda(t)\bigr)
 =u\bigl(\lambda(t)\bigr)-\log d-(\Re t)\bigl(L_d(\lambda)-\log d\bigr)
\end{gather*}
on $\bH_+$ is $\le L_d(\lambda(t))-L_d(\lambda(t))=0$ 
on $\{t\in\bC:(\Re t)(L_d(\lambda)-\log d)>\max\{0,\ell\}\}$
and vanishes if $\Re t\gg 1$, so by the maximum principle for subharmonic functions,
we have $u(\lambda)-L_d(\lambda)=0$ at this $\lambda\in E_d\cap\Omega_d(\ell)$.
Consequently,
\begin{gather*}
 u=L_d 
\end{gather*}
on $E_d\cap\Omega_d(\ell)$,
and in turn on $\Omega_d(\ell)$ by Lemma \ref{th:extend}.
Hence by the maximum principle for plurisubharmonic functions,
\begin{gather*}
u\le\lim_{0<\epsilon\searrow0}\sup_{\{L_d\le\log d+\max\{0,\ell\}+\epsilon\}}u
\le\lim_{0<\epsilon\searrow0}\sup_{\{L_d\le\log d+\max\{0,\ell\}+\epsilon\}}L_d
=\log d+\max\{0,\ell\}=L_d
\end{gather*}
on $\{L_d=\log d+\max\{0,\ell\}\}$, and in turn 
$u\le \log d+\max\{0,\ell\}$ on $\bC^{d-1}\setminus\Omega_d(\ell)$.

Let us conclude \eqref{eq:full}. If $\ell\le 0$, then
$u\le L_d$ on $C_d=\bC^{d-1}\setminus\Omega_d(\ell)$, and we are done by Lemma \ref{th:compare} and Theorem \ref{th:boundary} applied to this $u$, $v=L_d$, $\phi(\lambda)=\|\lambda\|^2$ on $\bC^{d-1}$,
and an open ball $\Omega$ containing $C_d$.

Alternatively if $\ell>0$, then 
for each $\lambda\in\{L_d<\log d+\ell\}=\bC^{d-1}\setminus\overline{\Omega_d(\ell)}$,
that is, $\log d\le L_d(\lambda)<\log d+\ell$,
letting $\lambda(t):\bH_+\to\bC^{d-1}$ be the holomorphic disk in $\bC^{d-1}$
passing through $\lambda(1)=\lambda$ obtained in Lemma \ref{th:stretch},
we have $\overline{\lambda((0,1])}\subset\{L_d\le L_d(\lambda)\}$
and $\lim_{r\searrow0}L_d(\lambda(r))=\log d$ by \eqref{eq:sliding}, 
so that $\bC^{d-1}\setminus\overline{\Omega_d(\ell)}$ is a domain in $\bC^{d-1}$
containing $C_d$.
Hence noting that $f_0(z)=z^d$ so $g_{f_0}=\log^+|\cdot|$,
$f_0^n(z)-z=z(z^{d^n-1}-1)$, $(f_0^n)'(z)=d^nz^{d^n-1}$, 
that $\lambda=0\in\Int C_d$, and that $|a_n|\ge e^{\ell n/2}\to\infty$ as $n\to\infty$, 
by Hartogs's lemma for plurisubharmonic functions, we compute
\begin{align*}
u(0)&\ge \lim_{j\to\infty}\frac{1}{\nu(n_j)}\log|p_{f,n_j}^*(0,w_{n_j})|
=\lim_{j\to\infty}\frac{1}{n_j}\frac{1}{\nu(n_j)}\cdot\nu(n_j)\log|d^{n_j}-w_{n_j}|\\
&=\lim_{j\to\infty}\Bigl(\log d+\frac{\log|a_{n_j}|}{n_j}+\frac{1}{n_j}
\log\Bigl|\frac{1}{a_{n_j}}-1\Bigr|\Bigr)=\log d+\ell,
\end{align*}
and we are done by the maximum principle for plurisubharmonic functions. 

Now the proof of Theorem \ref{th:main} is complete. \qed

\begin{proof}[Proof of Theorem \ref{th:unicrit}]
Under the holomorphic self-covering $\lambda=a^d$ of $\bC$, we work with the subfamily
$f_{0,a}(z)=z^d+a^d$ of the full critically marked family,
regarding all potential functions as subharmonic
functions on $\{c_0=\cdots=c_{d-2}=0\}\cong\bC$.
The intersections of $C_d$ and $E_d$ with $\{c_0=\cdots=c_{d-2}=0\}$
coincide with the inverse images of $M_d$ and $I_d$
under $(0,a)\mapsto a^d$, respectively,
so the extension step involving Lemma \ref{th:extend} is unnecessary.
We can replace Bassanelli--Berteloot's result for $w_n\equiv0$ with Levin's seminal work,
and Theorem \ref{th:boundary} with the boundary vanishing from Zakeri
and Graczyk--{\'S}wi\c{a}tek recalled in Subsection \ref{sec:unicrit}.
Recalling \eqref{eq:MP} and \eqref{eq:green},
we obtain \eqref{eq:mainuni} by an argument similar to that
in the proof of Theorem \ref{th:main}.
\end{proof}

\section{Proof of Theorem \ref{th:boundary}}\label{sec:boundary}

Let us denote by $\mathcal{P}_d^{\cm}$ the space of affine conjugacy classes of full critically
marked polynomials of degree $d$; see \cite[\S 5]{DF08}.
Let
$\pi_d:\bC^{d-1}\ni\lambda\mapsto[f_\lambda;0,c_1,\ldots,c_{d-2}]\in\mathcal{P}_d^{\cm}$
be the natural finite map (see \cite[Proposition 5.1]{DF08}).
For every $\lambda=(c_1,\ldots,c_{d-2},a)\in\bC^{d-1}$,
\begin{gather}\label{eq:fiber}
 \pi_d^{-1}(\pi_d(\lambda))
 =\bigl\{(\omega^d c_1,\ldots,\omega^d c_{d-2},\omega a):
 \omega\in\bC\quad\text{and}\quad\omega^{d(d-1)}=1\bigr\}.
\end{gather}
The corresponding conjugacy is $z\mapsto\omega^d z$.
Since $L_d$ is invariant under affine conjugacy, it descends to a continuous
plurisubharmonic function on $\mathcal{P}_d^{\cm}$, denoted by $\cL_d$ (see \cite[\S 5, proof of Proposition 5.1]{DF08}).
Let us define the bifurcation measure on $\mathcal{P}_d^{\cm}$ by
\begin{gather*}
 \mu_{\bif}:=(\rd\rd^c\cL_d)^{d-1},\quad\text{so that}\quad
 \pi_d^*\mu_{\bif}=\mu_d;
\end{gather*}
see \cite[Definition 6.10 and Proposition 6.14]{DF08} and \eqref{eq:MP}.
These measures do not charge analytic sets by \cite[Proposition 6.11]{DF08}.

\subsection{Combinatorial description of $\mu_{\bif}$}\label{sec:comb}

For every $\lambda\in C_d$ and every $\alpha\in\bR/\bZ$,
the external ray of angle $\alpha$ for $f_\lambda$ is
\begin{gather*}
 R_\alpha(\lambda):=
 \bigl\{\phi_{f_\lambda}^{-1}(e^{r+2\pi i\alpha}):r\in(0,+\infty)\bigr\}\subset\bC,
\end{gather*}
where $\phi_{f_\lambda}:\{g_{f_\lambda}>0\}\to\{|\zeta|>1\}$ is the
B\"ottcher coordinate from the proof of Lemma \ref{th:stretch},
tangent to the identity at $\infty$.
For every $\lambda\in C_d$ and every $\alpha\in\bQ/\bZ$,
$R_\alpha(\lambda)$ lands at
\begin{gather*}
 \ell_\alpha(\lambda):=
 \lim_{r\searrow0}\phi_{f_\lambda}^{-1}(e^{r+2\pi i\alpha})
 \in J_{f_\lambda},
\end{gather*}
which is called the landing point of $R_\alpha(\lambda)$.
The point $\ell_\alpha(\lambda)$ is preperiodic under $f_\lambda$, and the periodic cycle
it reaches is repelling or parabolic
(see Douady--Hubbard \cite[Expos\'e VIII, \S II.2, Proposition 2, p. 70]{DH84}).
For every $\lambda\in C_d$, the rational lamination of $f_\lambda$
is the relation (indeed, an equivalence relation)
\begin{gather*}
 \lambda_{\bQ}(f_\lambda):=
 \bigl\{(s,t)\in(\bQ/\bZ)^2:
 \ell_s(\lambda)=\ell_t(\lambda)\bigr\}
\end{gather*}
on $\bQ/\bZ$ (see \cite[\S 6.4, p. 92]{McMullen94}).
Let $\Cb$ be the space of critical portraits of degree $d$,
which is the totality of such elements
$\Theta=(\theta_0,\ldots,\theta_{d-2})\in(2^{\bR/\bZ})^{d-1}$
that
\begin{itemize}
\item $d\alpha=d\alpha'$ for every $j\in\{0,\ldots,d-2\}$ and any $\alpha,\alpha'\in\theta_j$,
\item either $\theta_i=\theta_j$ or $\theta_i\cap\theta_j=\emptyset$
for any distinct $i,j\in\{0,\ldots,d-2\}$,
\item for any distinct $i,j\in\{0,\ldots,d-2\}$,
if $\theta_i\cap\theta_j=\emptyset$, then
$\theta_j$ is contained in a component of
$(\bR/\bZ)\setminus\theta_i$,
\item $\#\{i\in\{0,\ldots,d-2\}:\theta_i=\theta_j\}=\#\theta_j-1$
for every $j\in\{0,\ldots,d-2\}$.
\end{itemize}
See \cite[\S 5 and Definition 7.7]{DF08}\footnotemark
\footnotetext{Our first three conditions correspond to the first, second, and fourth
conditions in \cite[Definition 7.7]{DF08}, respectively.
Our last condition implies the third condition there,
and also specifies the number of occurrences of each $\theta_j$.
The corresponding critical point has multiplicity $\#\theta_j-1$
and seems to need to be marked exactly $\#\theta_j-1$ times,
as required in \cite[\S 5]{DF08}.}.
For the itinerary equivalence relation
\begin{gather*}
 \Lambda_{\bQ}\bigl(\{\theta_0,\ldots,\theta_{d-2}\}\bigr)
\end{gather*}
associated to each $\Theta=(\theta_0,\ldots,\theta_{d-2})\in\Cb$,
see \cite[Definition 14.1]{Kiwi97}.

For every $\Theta=(\theta_0,\ldots,\theta_{d-2})\in\Cb$, the set
\begin{gather*}
 \Bigl\{\Cb\cap\prod_{j=0}^{d-2}2^{U_j}:
 \theta_j\subset U_j\subset\bR/\bZ
 \text{ and }U_j\text{ is open for every }j\in\{0,\ldots,d-2\}\Bigr\}
\end{gather*}
forms a neighborhood basis at $\Theta$ for the topology on $\Cb$;
see \cite[paragraph following Definition 7.7]{DF08}.
The metric space
\begin{gather*}
 \AP:=\bigl\{\{\alpha,\alpha'\}\subset\bR/\bZ:
 \alpha\ne\alpha'\quad\text{and}\quad d\alpha=d\alpha'\bigr\}
 \subset 2^{\bR/\bZ}
\end{gather*}
of angle pairs is equipped with the Hausdorff metric
induced by the usual distance on $\bR/\bZ$.
The subset $\AP^{d-1}\cap\Cb$ is open and dense in $\Cb$,
and its topology induced by $\Cb$ agrees with that induced by $\AP^{d-1}$.
Moreover, $\AP^{d-1}\cap\Cb$ is a nonempty open subset of $\AP^{d-1}$;
see \cite[Propositions 7.4 and 7.8]{DF08}.
For every $\{\alpha,\alpha'\}\in\AP$ and any $s,t\in(-1/(2d),1/(2d))$,
\begin{gather*}
 \{\alpha+s,\alpha'+t\}\in\AP
 \quad\Leftrightarrow\quad d(s-t)\in\bZ
 \quad\Leftrightarrow\quad s=t,
\end{gather*}
so that, for each $\Theta=(\theta_0,\ldots,\theta_{d-2})\in\AP^{d-1}$
and sufficiently small $\varepsilon>0$, the map
\begin{gather*}
 (-\varepsilon,\varepsilon)^{d-1}\ni(t_0,\ldots,t_{d-2})\mapsto
 (\theta_0+t_0,\ldots,\theta_{d-2}+t_{d-2})\in\AP^{d-1}
\end{gather*}
is a homeomorphism onto an open neighborhood of $\Theta$,
where $\theta_j+t_j:=\{\alpha+t_j:\alpha\in\theta_j\}$ for every
$j\in\{0,\ldots,d-2\}$.
The inverses of the displayed homeomorphisms are coordinate charts on $\AP^{d-1}$.
The Lebesgue measures $\rd t_0\cdots\rd t_{d-2}$
define a positive measure $m_{\AP}$ on $\AP^{d-1}$; indeed,
on the overlap of two such charts centered at this $\Theta$ and at
$\Theta'=(\theta'_0,\ldots,\theta'_{d-2})\in\AP^{d-1}$, their respective coordinates
$(t_0,\ldots,t_{d-2})$ and $(s_0,\ldots,s_{d-2})$ satisfy, for every
$j\in\{0,\ldots,d-2\}$,
\begin{gather*}
 \theta_j+t_j=\theta'_j+s_j
 \quad\Rightarrow\quad
 t_j-s_j\in\{u\in\bR:\theta_j+u=\theta'_j\},
\end{gather*}
and $\{u\in\bR:\theta_j+u=\theta'_j\}$ is discrete, so $t_j-s_j$ is locally constant and
\begin{gather*}
 \rd s_0\cdots\rd s_{d-2}=\rd t_0\cdots\rd t_{d-2}.
\end{gather*}
We compute
$m_{\AP}(\AP^{d-1})=((d-1)/2)^{d-1}<\infty$
(cf. \cite[Lemma 7.3 and \S 7.1.1]{DF08}).
Decomposing according to the two arcs determined by $\theta_0$, we compute
\begin{gather*}
 m_{\AP}(\AP^{d-1}\cap\Cb)
 =\frac{1}{2d^{d-2}}\sum_{k=1}^{d-1}
 \binom{d-2}{k-1}k^{k-2}(d-k)^{d-k-2}
 =\frac{1}{d}>0.
\end{gather*}
The exact values of these masses are not needed below.
Let us denote by $\mu_{\Cb}$ the trivial extension of
$m_{\AP}|_{\AP^{d-1}\cap\Cb}$ to $\Cb$, that is,
\begin{gather*}
 \mu_{\Cb}(A):=m_{\AP}(A\cap\AP^{d-1})
 \quad\text{for every Borel subset }A\text{ in }\Cb
\end{gather*}
(see \cite[Definitions 7.6 and 7.10]{DF08}).
For every $\Theta=(\theta_0,\ldots,\theta_{d-2})\in\Cb$, pick
$\lambda=\lambda_\Theta\in E_d$ such that, for every $j\in\{0,\ldots,d-2\}$,
\begin{gather*}
 g_{f_\lambda}(c_j(\lambda))=1\quad\text{and}\quad
 \theta_j=\Bigl\{\alpha\in\bR/\bZ:
 \lim_{s\searrow 1}\phi_{f_\lambda}^{-1}(e^{s+2\pi i\alpha})
 =c_j(\lambda)\Bigr\}.
\end{gather*}
Such a $\lambda$ exists by \cite[Proposition 3.8]{Goldberg94}, and
$\pi_d(\lambda)$ is uniquely determined by $\Theta$
(see \cite[Lemma 9.7]{Kiwi97}).
The holomorphic disk $\lambda_\Theta(t):\bH_+\to\bC^{d-1}$ constructed in the proof of
Lemma \ref{th:stretch}, passing through $\lambda_\Theta(1)=\lambda_\Theta$,
preserves the external angles for $t=r>0$ and satisfies
$g_{f_{\lambda_\Theta(r)}}(c_j(\lambda_\Theta(r)))=r$
for every $r>0$ and every $j\in\{0,\ldots,d-2\}$
(see also \cite[\S 8, Proposition 8.4]{BrannerHubbard88} and
\cite[Remark 7.14]{DF08}).
For every $r>0$, the class $\pi_d(\lambda_\Theta(r))$ is independent of the choice of $\lambda_\Theta$.

\subsection{Proof of Theorem \ref{th:boundary}}

Let $U$ be a component of $\Int C_d$.
All marked critical points $c_0,\ldots,c_{d-2}$ of $f$ are passive on $U$, that is,
$\{\lambda\mapsto f_\lambda^n(c_j(\lambda)):U\to\bP^1\}_{n\ge 0}$ is a normal family
for every $0\le j\le d-2$,
since these iterates are locally uniformly bounded in $\bC$ on $U$, and then
$(f_\lambda)_{\lambda\in U}$ is $J$-stable (see \cite[Lemma III.2 and Theorem B]{MSS}).
If $f_\lambda$ had an indifferent cycle of formal exact period $n\in\bN$ having multiplier $w\in\partial\bD$ for some
$\lambda\in U$, that would persist locally with the same multiplier
(see \cite[\S I.6 and Theorem B]{MSS}),
so $p_{f,n}^*(\cdot,w)\equiv 0$ on $\bC^{d-1}$, which is impossible since
$p_{f,n}^*(0,w)\ne 0$ for this $w$
(see the computation at the end of Section \ref{sec:complex}).
Thus the number of attracting cycles of $f_\lambda$ in $\bC$ of each (formal) exact period
is constant for $\lambda\in U$.

If $f_\lambda$ has an attracting cycle in $\bC$ of (formal) exact period $n$ for some $\lambda\in U$,
then $\partial U\subset\bigcup_{|w|\le 1}
\bigl\{\lambda\in\bC^{d-1}:p_{f,n}^*(\lambda,w)=0\bigr\}$.
By \cite[Corollary 11]{DF08}, $\mu_{\bif}$-almost every polynomial
has only repelling cycles in $\bC$, so
$\mu_{\bif}(\pi_d(\partial U))=\mu_d(\partial U)=0$.

It remains the case that all cycles of $f_\lambda$ in $\bC$ are repelling for every
$\lambda\in U$. In this case, $c_j(\lambda)\in J_{f_\lambda}$ for every
$\lambda\in U$ and every $0\le j\le d-2$ (see \cite[Theorem 1, p. 404]{Sullivan85}).
Let us see that $c_j(\lambda)$ is not preperiodic under $f_\lambda$ for every
$\lambda\in U$ and every $0\le j\le d-2$; for, otherwise,
$f_{\lambda_0}^n(c_j(\lambda_0))=f_{\lambda_0}^m(c_j(\lambda_0))$
for some $j\in\{0,\ldots,d-2\}$, some $\lambda_0\in U$, and some $n,m\in\bN\cup\{0\}$, $n>m$, and then under the above standing assumption,
$C:=\bigl|(f_{\lambda_0}^{n-m})'(f_{\lambda_0}^m(c_j(\lambda_0)))\bigr|>1$.
Hence if $|\lambda-\lambda_0|\ll 1$, then
$|z-z'|\le(2/(C+1))|f_\lambda^{n-m}(z)-f_\lambda^{n-m}(z')|$
for any $z,z'\in\bC$ close enough to $f_{\lambda_0}^m(c_j(\lambda_0))$, 
and then by the normality of $\{\lambda\mapsto f_\lambda^n(c_j(\lambda)):U\to\bP^1\}_{n\ge 0}$ at $\lambda_0$,
\begin{gather*}
\bigl|f_\lambda^n\bigl(c_j(\lambda)\bigr)-f_\lambda^m\bigl(c_j(\lambda)\bigr)\bigr|
\le \Bigl(\frac{2}{C+1}\Bigr)^k
 \sup_{\ell\ge 0}\bigl|f_\lambda^{(\ell+1)(n-m)}\bigl(f_\lambda^m(c_j(\lambda))\bigr)-f_\lambda^{\ell(n-m)}\bigl(f_\lambda^m(c_j(\lambda))\bigr)\bigr|\to 0
\end{gather*}
as $k\to\infty$ (cf.\ \cite[Lemma 1.3 and the proof]{DF08}).
By the identity theorem, we must have
$f_\lambda^n(c_j(\lambda))\equiv f_\lambda^m(c_j(\lambda))$
on $\bC^{d-1}\supset E_d$, which is impossible.

Consider the following two conditions on each $\lambda\in C_d$:
\begin{enumerate}
\renewcommand{\theenumi}{\roman{enumi}}
\renewcommand{\labelenumi}{(\theenumi)}
\item\label{cond:critical}
All marked critical points $c_0,\ldots,c_{d-2}$ of $f_\lambda$
are not preperiodic under $f_\lambda$.
\item\label{cond:cycles}
All cycles of $f_\lambda$ in $\bC$ are repelling.
\end{enumerate}
Pick any $\lambda^*\in\partial U$ such that 
(\ref{cond:critical}) and (\ref{cond:cycles}) hold.
Let us first see that for every $\lambda'\in U$,
\begin{gather}
 \lambda_{\bQ}(f_{\lambda^*})=\lambda_{\bQ}(f_{\lambda'});
\end{gather}
indeed, for every $\alpha\in\bQ/\bZ$, take any $n,m\in\bN\cup\{0\}$, $n>m$, such that
$d^n\alpha=d^m\alpha$.
By \cite[Lemma 17.1]{Kiwi97}, for every
$\lambda_0\in U\cup\{\lambda^*\}$, the unique local holomorphic 
function $z(\lambda)$ solving the functional equation
$f_\lambda^n=f_\lambda^m$ near $\lambda_0$ such that
$z(\lambda_0)=\ell_\alpha(\lambda_0)$ 
coincides with $\ell_\alpha(\lambda)$ near $\lambda_0$,
so $\ell_\alpha$ is holomorphic on $U$ and extends holomorphically
to some open neighborhood of $\lambda^*$ in $\bC^{d-1}$. 
Hence for any distinct $s,t\in\bQ/\bZ$, choosing common $n,m\in\bN\cup\{0\}$, $n>m$,
such that $d^n s=d^m s$ and $d^n t=d^m t$, the set
$\bigl\{\lambda\in U\cup\{\lambda^*\}:
 \ell_s(\lambda)=\ell_t(\lambda)\bigr\}$
is clopen in $U\cup\{\lambda^*\}$ by the above common description of $\ell_s$ 
and $\ell_t$, so for every $\lambda'\in U$, since
$U\cup\{\lambda^*\}$ is connected, we have
$\ell_s(\lambda^*)=\ell_t(\lambda^*)
\Leftrightarrow
 \ell_s(\lambda')=\ell_t(\lambda')$, and we are done.

Using the holomorphic disk $\lambda_\Theta(t)$ chosen above for each $\Theta\in\Cb$
in Subsection \ref{sec:comb}, set
\begin{gather*}
 \Cb_1:=\Biggl\{\Theta\in\Cb:
\begin{gathered}
 \lim_{r\searrow0}\pi_d(\lambda_\Theta(r))=\pi_d(\lambda)\text{ for some }\lambda\in C_d
 \text{ such that (\ref{cond:critical}) and (\ref{cond:cycles}) hold}\\
 \text{and that all marked critical points of }f_{\lambda}\text{ are simple}
\end{gathered}
\Biggr\}.
\end{gather*}
By \cite[Proposition 7.19 and \S 7.6, p. 1027]{DF08},
$\Cb_1$ is $\mu_{\Cb}$-measurable and $\mu_{\Cb}(\Cb\setminus\Cb_1)=0$.
The map
\begin{gather*}
 \Phi_1:\Cb_1\ni\Theta
 \mapsto
 \lim_{r\searrow0}\pi_d\bigl(\lambda_\Theta(r)\bigr)
 \in\mathcal{P}_d^{\cm}
\end{gather*}
is $\mu_{\Cb}|\Cb_1$-measurable (see \cite[\S 7.3, proof of Theorem 8]{DF08}).
Then $\mu_{\bif}=(d-2)!(\Phi_1)_*(\mu_{\Cb}|\Cb_1)$
(\cite[Theorem 8]{DF08}).

Let us next see that
\begin{gather}
 \#\biggl(
 \bigcup_{\lambda^*\in\partial U}
 \Bigl\{\Theta\in\Cb_1:
 \lim_{r\searrow0}\pi_d(\lambda_\Theta(r))
 =\pi_d(\lambda^*)\Bigr\}
 \biggr)<\infty;\label{eq:finiteness}
\end{gather}
indeed, for every $\lambda^*\in\partial U$ 
and every $\Theta=(\theta_0,\ldots,\theta_{d-2})\in\Cb_1$
such that
$\lim_{r\searrow0}\pi_d(\lambda_{\Theta}(r))=\pi_d(\lambda^*)$,
by the properness of $\pi_d$, the accumulation set
of $\lambda_{\Theta}(r)$ as $r\searrow0$ is a nonempty compact subset of the finite fiber
$\pi_d^{-1}(\pi_d(\lambda^*))$ and in turn
$\lambda_{\Theta}(r)$ tends to a point in $\pi_d^{-1}(\pi_d(\lambda^*))$
as $r\searrow0$ by the continuity of $r\mapsto\lambda_{\Theta}(r)$.
Hence by \eqref{eq:fiber}, for some $k\in\{0,\ldots,d-2\}$,
setting $R_k(z)=e^{2\pi i k/(d-1)}z$,
$f_{\lambda_{\Theta}(r)}
 \to R_k\circ f_{\lambda^*}\circ R_k^{-1}$
as $r\searrow0$ in the space of degree $d$ complex polynomials,
and moreover, by \cite[Definition 7.6 and Theorem 14.4]{Kiwi97}, 
for any $\lambda'\in U$, we have
\begin{gather*}
 \Lambda_{\bQ}(\{\theta_0,\ldots,\theta_{d-2}\})=\lambda_{\bQ}(R_k\circ f_{\lambda^*}\circ R_k^{-1})=\lambda_{\bQ}(R_k\circ f_{\lambda'}\circ R_k^{-1}).
\end{gather*}
Hence fixing $\lambda'\in U$,
\begin{multline*}
 \bigcup_{\lambda^*\in\partial U}
 \Bigl\{\Theta\in\Cb_1:
 \lim_{r\searrow0}\pi_d\bigl(\lambda_\Theta(r)\bigr)
 =\pi_d(\lambda^*)\Bigr\}\\
 \subset\bigcup_{k=0}^{d-2}
 \biggl\{(\theta_0,\ldots,\theta_{d-2})\in\Cb:
 \Lambda_{\bQ}(\{\theta_0,\ldots,\theta_{d-2}\})
 =\lambda_{\bQ}(R_k\circ f_{\lambda'}\circ R_k^{-1})
 \biggr\},
\end{multline*}
for which 
recalling that each distinct component
$\theta$ of an element of $\Cb$ occurs $\#\theta-1$ times
among all the $d-1$ components, by \cite[Theorem 14.3]{Kiwi97}, we have
\begin{gather*}
 \#\text{RHS}\le (d-1)!\sum_{k=0}^{d-2}
 \#\biggl\{
 \{\theta_0,\ldots,\theta_{d-2}\}:
 \begin{gathered}
 (\theta_0,\ldots,\theta_{d-2})\in\Cb\quad \text{and}\\
 \Lambda_{\bQ}(\{\theta_0,\ldots,\theta_{d-2}\})
 =\lambda_{\bQ}(R_k\circ f_{\lambda'}\circ R_k^{-1})
 \end{gathered}
 \biggr\}<\infty,
\end{gather*}
and we are done.
Once \eqref{eq:finiteness} is at our disposal,
noting that $\mu_{\Cb}$ has no atoms by \cite[Definitions 7.6 and 7.10]{DF08},
we have
\begin{multline*}
 \mu_d(\partial U)
=\mu_{\bif}\bigl(\pi_d(\partial U)\bigr)
=(d-2)!\bigl((\Phi_1)_*(\mu_{\Cb}|\Cb_1)\bigr)\bigl(\pi_d(\partial U)\bigr)\\
=(d-2)!\mu_{\Cb}\biggl(
 \bigcup_{\lambda^*\in\partial U}
 \Bigl\{\Theta\in\Cb_1:
 \lim_{r\searrow0}\pi_d\bigl(\lambda_\Theta(r)\bigr)
 =\pi_d(\lambda^*)\Bigr\}
 \biggr)=0.
\end{multline*}
Now the proof of Theorem \ref{th:boundary} is complete.\qed

\appendix
\section{Proofs of the thermodynamic preliminaries}\label{th:thermo}

For more details on the thermodynamic formalism used below,
see \cite[Chapter 5 and \S 6.4]{PU10}. We show the assertions in Subsection \ref{sec:thermo}
for the full critically marked polynomial family $f$ of degree $d$ in Subsection \ref{sec:crit}; the proof for the unicritical polynomial family in Subsection \ref{sec:unicrit} is similar. Fix $\lambda_0\in E_d$, and
we use the $d$-symbol coding of the expanding Cantor Julia set and its
continuation by holomorphic motion from $\lambda$-lemma. Pick any $\theta\in(0,1)$. Then
the metric space $\Sigma_d:=\{0,\ldots,d-1\}^{\{0\}\sqcup\bN}$ equipped with 
the metric 
\begin{gather*}
 d_\theta((\omega_n),(\omega_n'))
=\theta^{\min\{j\ge 0:\omega_j\ne\omega_j'\}} 
\end{gather*}
under the convention $\min\emptyset=+\infty$ is compact, 
and the shift selfmap 
\begin{gather*}
 \sigma((\omega_0,\omega_1,\ldots))=(\omega_1,\omega_2,\ldots)
\end{gather*}
of $\Sigma_d$ is $\theta^{-1}$-Lipschitz continuous.

\begin{notation}
Let us write $\overline{k_0,\ldots,k_{n-1}}=((k_0,\ldots,k_{n-1}),\ldots)\in\Sigma_d$ 
fixed by $\sigma^n$ for each $n\in\bN$ and each $(k_0,k_1,\ldots,k_{n-1})\in\{0,\ldots,d-1\}^n$, and set $i\omega:=(i,\omega_0,\omega_1,\ldots)\in\Sigma_d$ for
each $i\in\{0,\ldots,d-1\}$ and each $\omega=(\omega_0,\omega_1,\ldots)\in\Sigma_d$.
\end{notation}

For every open disk $U_0=\{\lambda\in\bC^{d-1}:\|\lambda-\lambda_0\|<r_0\}\Subset E_d$
for $0<r_0\ll 1$,
there are homeomorphisms $\psi_\lambda:\Sigma_d\to J_{f_\lambda}$,
$\lambda\in\overline{U_0}$, such that
\begin{gather*}
 f_\lambda\circ\psi_\lambda=\psi_\lambda\circ\sigma\quad\text{on }\Sigma_d.
\end{gather*}
The map $\psi_{\lambda}(\omega)$ is jointly continuous on 
$(\lambda,\omega)\in\overline{U_0}\times\Sigma_d$, and
for every $\omega\in\Sigma_d$, $\lambda\mapsto\psi_\lambda(\omega)$ is holomorphic on $U_0$. 
By the expanding property of $f_\lambda$ on $J_{f_\lambda}$ locally uniformly
on $E_d$, diminishing $U_0$ if necessary, we also have $m\in\bN$ such that
\begin{gather*}
\inf_{\lambda\in U_0}\inf_{\omega\in\Sigma_d}|(f_\lambda^m)'(\psi_\lambda(\omega))|>1.
\end{gather*}
Using the Taylor theorem,
if $N\gg 1$, then for any $\omega,\omega'\in\Sigma_d$ such that 
$\min\{j\in\bN\sqcup\{0\}:\omega_j\neq\omega_j'\}\ge N$, we have
\begin{gather*}
 |\psi_\lambda(\sigma^m\omega)-\psi_\lambda(\sigma^m\omega')|
 \ge \Bigl(\inf_{\lambda\in U_0}\inf_{\omega\in\Sigma_d}
 |(f_\lambda^m)'(\psi_\lambda(\omega))|\Bigr)^{1/2}\cdot
 |\psi_\lambda(\omega)-\psi_\lambda(\omega')|,\quad\text{and in turn}\\
 |\psi_\lambda(\omega)-\psi_\lambda(\omega')|
 \le \Bigl(\inf_{\lambda\in U_0}\inf_{\omega\in\Sigma_d}
 |(f_\lambda^m)'(\psi_\lambda(\omega))|\Bigr)^{-\lfloor(\min\{j\in\bN\sqcup\{0\}:\omega_j\neq\omega_j'\}-N)/m\rfloor/2}
 \sup_{\lambda\in U_0}\operatorname{diam}J_{f_\lambda}
\end{gather*}
inductively.
In particular,
\begin{multline*}
 |\psi_\lambda(\omega)-\psi_\lambda(\omega')|\le \Bigl(\inf_{\lambda\in U_0}\inf_{\omega\in\Sigma_d}
 |(f_\lambda^m)'(\psi_\lambda(\omega))|\Bigr)^{-(\min\{j\in\bN\sqcup\{0\}:\omega_j\neq\omega_j'\}-N-m)/(2m)}
 \sup_{\lambda\in U_0}\operatorname{diam}J_{f_\lambda}\\
 =\Bigl(\inf_{\lambda\in U_0}\inf_{\omega\in\Sigma_d}
 |(f_\lambda^m)'(\psi_\lambda(\omega))|\Bigr)^{(N+m)/(2m)}
 \Bigl(\sup_{\lambda\in U_0}\operatorname{diam}J_{f_\lambda}\Bigr)\\
 \times\theta^{\bigl(-(\min\{j\in\bN\sqcup\{0\}:\omega_j\neq\omega_j'\})
 \log_\theta\inf_{\lambda\in U_0}\inf_{\omega\in\Sigma_d}
 |(f_\lambda^m)'(\psi_\lambda(\omega))|\bigr)/(2m)}\\
\le \Bigl(\inf_{\lambda\in U_0}\inf_{\omega\in\Sigma_d}
 |(f_\lambda^m)'(\psi_\lambda(\omega))|\Bigr)^{(N+m)/(2m)}
 \Bigl(\sup_{\lambda\in U_0}\operatorname{diam}J_{f_\lambda}\Bigr)
 \times d_\theta(\omega,\omega')^\beta,
\end{multline*}
where
\begin{gather*}
 \beta=\beta_\theta:=\min\Bigl\{1,
 \frac{\log\inf_{\lambda\in U_0}\inf_{\omega\in\Sigma_d}
 |(f_\lambda^m)'(\psi_\lambda(\omega))|}{-2m\log\theta}\Bigr\}>0,
\end{gather*}
and on the other hand, for $N\gg 1$, 
if $\min\{j\in\bN\sqcup\{0\}:\omega_j\neq\omega_j'\}<N$, then
we still have
\begin{multline*}
 |\psi_\lambda(\omega)-\psi_\lambda(\omega')|
 \le \Bigl(\inf_{\lambda\in U_0}\inf_{\omega\in\Sigma_d}
 |(f_\lambda^m)'(\psi_\lambda(\omega))|\Bigr)^{N/(2m)}
 \Bigl(\sup_{\lambda\in U_0}\operatorname{diam}J_{f_\lambda}\Bigr)\\
 \times\theta^{\bigl(-\min\{j\in\bN\sqcup\{0\}:\omega_j\ne\omega_j'\}
 \log_\theta\inf_{\lambda\in U_0}\inf_{\omega\in\Sigma_d}
 |(f_\lambda^m)'(\psi_\lambda(\omega))|\bigr)/(2m)}\\
 \le \Bigl(\inf_{\lambda\in U_0}\inf_{\omega\in\Sigma_d}
 |(f_\lambda^m)'(\psi_\lambda(\omega))|\Bigr)^{(N+m)/(2m)}
 \Bigl(\sup_{\lambda\in U_0}\operatorname{diam}J_{f_\lambda}\Bigr)
 \times d_\theta(\omega,\omega')^\beta.
\end{multline*}
Namely, $\psi_\lambda(\omega)$ is H\"older continuous on $\Sigma_d$
uniformly on $\lambda\in U_0$. Then 
we compute not only
\begin{align*}
 &\bigl\|\log|f_\lambda'\circ\psi_\lambda|\bigr\|_{C^\beta(\Sigma_d)}
 =\bigl\|\log|f_\lambda'\circ\psi_\lambda|\bigr\|_{C^0(\Sigma_d)}
 +\sup_{\omega\ne\omega'}
 \frac{\bigl|\log|f_\lambda'(\psi_\lambda(\omega))|
 -\log|f_\lambda'(\psi_\lambda(\omega'))|\bigr|}
 {d_\theta(\omega,\omega')^\beta}\\
 \le& \bigl\|\log|f_\lambda'\circ\psi_\lambda|\bigr\|_{C^0(\Sigma_d)}
 +
 \sup_{\omega\ne\omega'}
 \frac{1}{d_\theta(\omega,\omega')^\beta}
 \frac{|f_\lambda'(\psi_\lambda(\omega))
 -f_\lambda'(\psi_\lambda(\omega'))|}{\inf_{z\in J_{f_\lambda}}|f_\lambda'(z)|}\\
 \le& \sup_{\lambda\in U_0}
 \bigl\|\log|f_\lambda'\circ\psi_\lambda|\bigr\|_{C^0(\Sigma_d)}\\
 &+\frac{\sup_{\lambda\in U_0}
 \sup_{z,z'\in J_{f_\lambda}:\,\text{distinct}}
 |f_\lambda'(z)-f_\lambda'(z')|/|z-z'|}
 {\inf_{\lambda\in U_0}\inf_{z\in J_{f_\lambda}}|f_\lambda'(z)|}
 \sup_{\lambda\in U_0}\sup_{\omega\ne\omega'}
 \frac{|\psi_\lambda(\omega)-\psi_\lambda(\omega')|}
 {d_\theta(\omega,\omega')^\beta}
 <\infty
\end{align*}
for every $\lambda\in U_0$, so that
\begin{gather}
 C:=\sup_{\lambda\in U_0}
 \bigl\|\log|f_\lambda'\circ\psi_\lambda|\bigr\|_{C^\beta(\Sigma_d)}
 <\infty,\label{eq:bound}
\end{gather}
but also for every $\alpha\in(0,\beta)$, as $\lambda\to\lambda_0$,
\begin{align*}
 &\bigl\|\log|f_\lambda'\circ\psi_\lambda|
 -\log|f_{\lambda_0}'\circ\psi_{\lambda_0}|\bigr\|_{C^\alpha(\Sigma_d)}
 =\bigl\|\log|f_\lambda'\circ\psi_\lambda|
 -\log|f_{\lambda_0}'\circ\psi_{\lambda_0}|\bigr\|_{C^0(\Sigma_d)}\\
 &+\sup_{\omega\ne\omega'}\frac{1}{d_\theta(\omega,\omega')^\alpha}
 \Bigl|\bigl(\log|f_\lambda'\circ\psi_\lambda|
 -\log|f_{\lambda_0}'\circ\psi_{\lambda_0}|\bigr)(\omega)
 -\bigl(\log|f_\lambda'\circ\psi_\lambda|
 -\log|f_{\lambda_0}'\circ\psi_{\lambda_0}|\bigr)(\omega')\Bigr|\\
  &\le \bigl\|\log|f_\lambda'\circ\psi_\lambda|
 -\log|f_{\lambda_0}'\circ\psi_{\lambda_0}|\bigr\|_{C^0(\Sigma_d)}\\
 &\quad+\sup_{\omega\ne\omega'}\frac{1}{d_\theta(\omega,\omega')^\alpha}
 \Bigl(2\bigl\|\log|f_\lambda'\circ\psi_\lambda|
 -\log|f_{\lambda_0}'\circ\psi_{\lambda_0}|\bigr\|_{C^0(\Sigma_d)}\Bigr)^{1-\alpha/\beta}
 \bigl(2C d_\theta(\omega,\omega')^\beta\bigr)^{\alpha/\beta}\\
 &\le(1+2C^{\alpha/\beta})\bigl\|\log|f_\lambda'\circ\psi_\lambda|
 -\log|f_{\lambda_0}'\circ\psi_{\lambda_0}|\bigr\|_{C^0(\Sigma_d)}^{1-\alpha/\beta}\\
&\le(1+2C^{\alpha/\beta})
\biggl(\frac{\| |f_\lambda'\circ\psi_\lambda|-|f_{\lambda_0}'\circ\psi_{\lambda_0}|\|_{C^0(\Sigma_d)}}{\inf_{\lambda\in U_0}\inf_{z\in J_{f_\lambda}}|f_\lambda'(z)|}\biggr)^{1-\alpha/\beta}\to 0,
\end{align*}
that is, $U_0\ni\lambda\mapsto\log|f_\lambda'\circ\psi_\lambda|\in C^\alpha(\Sigma_d)$ 
is continuous at $\lambda=\lambda_0$.

For every $\lambda\in U_0$ and every $t\in\bC$, 
the transfer operator $\cL_{\lambda,t}\in L(C^0(J_{f_\lambda}))$
is defined as in \eqref{eq:ruelle}. For each $\lambda\in U_0$, 
the pullback 
$\psi_\lambda^*\in L(C^0(J_{f_\lambda}),C^0(\Sigma_d))$ and the pushforward
$(\psi_\lambda)_*\in L(C^0(\Sigma_d),C^0(J_{f_\lambda}))$
are mutually inverse isometries, and 
$\psi_\lambda^*\circ\cL_{\lambda,t}\circ(\psi_\lambda)_*\in L(C^0(\Sigma_d))$
is the transfer operator with the weight $-t\log|f_\lambda'\circ\psi_{\lambda}|$.
 
For every $(\lambda,t)\in U_0\times\bC$,
every $h\in C^\alpha(\Sigma_d)$, and every $k\in\{0,1\}$, we have the function
\begin{gather*}
 \bigl((\psi_\lambda^*\circ(\partial_t^k\cL_{\lambda,t})\circ(\psi_\lambda)_*)h\bigr)(\omega)
 =\sum_{i=0}^{d-1}
 \bigl(-\log|f_\lambda'(\psi_\lambda(i\omega))|\bigr)^k
 e^{-t\log|f_\lambda'(\psi_\lambda(i\omega))|}h(i\omega)\in C^0(\Sigma_d),
\end{gather*}
and indeed, we have the map 
\begin{gather*}
 U_0\times\bC\ni(\lambda,t)\mapsto\psi_\lambda^*\circ(\partial_t^k\cL_{\lambda,t})\circ(\psi_\lambda)_*\in L(C^\alpha(\Sigma_d)), 
\end{gather*}
which is not only jointly continuous at $(\lambda_0,t_0)$ for every $t_0\in\bC$,
that is, 
\begin{gather}\label{eq:opcont}
 \bigl\|\psi_\lambda^*\circ(\partial_t^k\cL_{\lambda,t})\circ(\psi_\lambda)_*
 -\psi_{\lambda_0}^*\circ(\partial_t^k\cL_{\lambda_0,t_0})\circ(\psi_{\lambda_0})_*
 \bigr\|_{L(C^\alpha(\Sigma_d))}\to 0
 \quad\text{as }(\lambda,t)\to(\lambda_0,t_0),
\end{gather}
but also holomorphic in $t\in\bC$ for every fixed $\lambda\in U_0$;
for, when $k=0$, by $d_\theta(i\omega,i\omega')=\theta d_\theta(\omega,\omega')$,
the submultiplicativity of $\|\cdot\|_{C^\alpha(\Sigma_d)}$
and \eqref{eq:bound}, we compute
\begin{align*}
 &\bigl\|(\psi_\lambda^*\circ\cL_{\lambda,t}\circ(\psi_\lambda)_*)h
 \bigr\|_{C^\alpha(\Sigma_d)}
 =\Bigl\|\sum_{i=0}^{d-1}
 \bigl(e^{-t\log|f_\lambda'\circ\psi_\lambda|}h\bigr)(i\,\cdot)
 \Bigr\|_{C^\alpha(\Sigma_d)}\\
 &\le \sum_{i=0}^{d-1}\bigl\|e^{-t\log|f_\lambda'\circ\psi_\lambda|}h
 \bigr\|_{C^\alpha(\Sigma_d)}
 \le d\bigl\|e^{-t\log|f_\lambda'\circ\psi_\lambda|}
 \bigr\|_{C^\alpha(\Sigma_d)}\|h\|_{C^\alpha(\Sigma_d)}\\
 &\le d e^{|t|\|\log|f_\lambda'\circ\psi_\lambda|\|_{C^\alpha(\Sigma_d)}}
 \cdot\|h\|_{C^\alpha(\Sigma_d)}
 \le d e^{|t|C}\cdot\|h\|_{C^\alpha(\Sigma_d)},
\end{align*}
and the continuity of $\log|f_\lambda'\circ\psi_\lambda|$
at $\lambda=\lambda_0$ in $C^\alpha(\Sigma_d)$ above
yields the regularity assertions for $k=0$. The case of $k=1$ follows from 
\begin{gather*}
 \bigl(\psi_\lambda^*\circ(\partial_t\cL_{\lambda,t})\circ(\psi_\lambda)_*\bigr)h
 =- \bigl(\psi_\lambda^*\circ\cL_{\lambda,t}\circ(\psi_\lambda)_*\bigr)\bigl((\log|f_\lambda'\circ\psi_\lambda|)h\bigr)
\end{gather*}
and the submultiplicativity of $\|\cdot\|_{C^\alpha(\Sigma_d)}$.

For each $\lambda\in U_0$ and each $t\in\bR$, 
by \eqref{eq:bound} and the Perron--Frobenius--Ruelle theorem 
(see e.g.\ \cite[Theorem 4.1]{Zins00}) applied to the transfer operator
$\psi_\lambda^*\circ\cL_{\lambda,t}\circ(\psi_\lambda)_*\in L(C^\alpha(\Sigma_d))$, 
there are a simple eigenvalue $\Lambda_{\lambda,t}$ of 
$\psi_\lambda^*\circ\cL_{\lambda,t}\circ(\psi_\lambda)_*$ and the unique strictly positive  eigenfunction $h^\Sigma_{\lambda,t}\in C^\alpha(\Sigma_d)$ and the unique probability eigenmeasure $m^\Sigma_{\lambda,t}$ on $\Sigma_d$ of $\psi_\lambda^*\circ\cL_{\lambda,t}\circ(\psi_\lambda)_*$ and ${}^t(\psi_\lambda^*\circ\cL_{\lambda,t}\circ(\psi_\lambda)_*)$, respectively, belonging to $\Lambda_{\lambda,t}$ and
normalized as $m^\Sigma_{\lambda,t}(1)=m^\Sigma_{\lambda,t}(h^\Sigma_{\lambda,t})=1$;
using the positivity of $h^\Sigma_{\lambda,t}$, such $\Lambda_{\lambda,t}$ is unique, and
moreover (see e.g.\ \cite[Theorem 4.5]{Zins00}), we have the spectral gap or the aperiodicity
\begin{gather*}
 \Lambda_{\lambda,t}>0\quad\text{and}\quad
\sup\bigl\{|z|:
 z\in\spe_{C^\alpha(\Sigma_d)}
 (\psi_\lambda^*\circ\cL_{\lambda,t}\circ(\psi_\lambda)_*)
 \setminus\{\Lambda_{\lambda,t}\}\bigr\}
 <\Lambda_{\lambda,t}.
\end{gather*}
Then for every $t_0\in\bR$,
by \eqref{eq:opcont} and 
the Kato--Rellich perturbation theorem (see e.g.\ \cite[Theorem 4.24]{Zins00})
applied to
$\psi_{\lambda_0}^*\circ\cL_{\lambda_0,t_0}\circ(\psi_{\lambda_0})_*
\in L(C^\alpha(\Sigma_d))$, the map
\begin{gather*}
 U_0\times\bR\ni(\lambda,t)\mapsto
 \bigl(\Lambda_{\lambda,t},h^\Sigma_{\lambda,t}\bigr)
 \in\bC\times C^\alpha(\Sigma_d)
\end{gather*}
extends jointly continuously in $(\lambda,t)$ and holomorphically in $t$
near $(\lambda_0,t_0)$, so in particular we have
the joint continuity of $\partial_t\log\Lambda_{\lambda,t}$ there, as stated in Subsection \ref{sec:thermo}. Consequently, 
for $\lambda\in U_0$ and every $t\in\bR$,
by the Perron--Frobenius--Ruelle theorem again, 
\begin{gather*}
 h_{\lambda,t}=(\psi_\lambda)_*h^\Sigma_{\lambda,t}\quad\text{and}\quad
m_{\lambda,t}=(\psi_\lambda)_*m^\Sigma_{\lambda,t}.
\end{gather*}
Let us see the remaining assertions.
Pick any $\lambda\in U_0$. For every $t_0\in\bR$, 
differentiating both sides of the
eigenfunction equation 
$\Lambda_{\lambda,t}h^\Sigma_{\lambda,t}=(\psi_\lambda^*\circ(\cL_{\lambda,t})
 \circ(\psi_\lambda)_*)h^\Sigma_{\lambda,t}$ in $t$ near $t=t_0$, we have
\begin{align*}
(\partial_t\Lambda_{\lambda,t})h^\Sigma_{\lambda,t}
&=
 \bigl(\psi_\lambda^*\circ(\partial_t\cL_{\lambda,t})
 \circ(\psi_\lambda)_*\bigr)h^\Sigma_{\lambda,t}
 +\bigl(\psi_\lambda^*\circ\cL_{\lambda,t}
 \circ(\psi_\lambda)_*\bigr)\partial_t h^\Sigma_{\lambda,t}
-\Lambda_{\lambda,t}\partial_t h^\Sigma_{\lambda,t},
\end{align*}
and integrating both sides at $t=t_0$
against $m^\Sigma_{\lambda,t_0}$,
\begin{multline*}
 (\partial_t\Lambda_{\lambda,t})\bigm|_{t=t_0}\cdot 1
=m^\Sigma_{\lambda,t_0}\bigl(
 \bigl(\psi_\lambda^*\circ(\partial_t\cL_{\lambda,t}\big|_{t=t_0})
\circ(\psi_\lambda)_*\bigr)h^\Sigma_{\lambda,t_0}\bigr)+0\\
 =m^\Sigma_{\lambda,t_0}\Bigl(
 \bigl(\psi_\lambda^*\circ\cL_{\lambda,t_0}\circ(\psi_\lambda)_*\bigr)
 \bigl(-\log|f_\lambda'\circ\psi_\lambda|\cdot h^\Sigma_{\lambda,t_0}\bigr)\Bigr)
 =\Lambda_{\lambda,t_0}\cdot
 \int_{\Sigma_d}-\log|f_\lambda'\circ\psi_\lambda|
 h^\Sigma_{\lambda,t_0}\,m^\Sigma_{\lambda,t_0},
\end{multline*} 
so in particular, when $t_0=0$, we have the desired equality \eqref{eq:deriv}.
Finally, for every $\lambda\in U_0$ and
any distinct $\omega,\omega'\in\Sigma_d$, 
setting $\min\{j\in\bN\sqcup\{0\}:\omega_j\neq\omega_j'\}=:n$,
by \eqref{eq:bound} and the chain rule, we have
\begin{gather*}
 \biggl|\log\frac{|(f_\lambda^n)'(\psi_\lambda(\omega))|}{|(f_\lambda^n)'(\psi_\lambda(\omega'))|}\biggr|
 \le \sum_{j=0}^{n-1}
 \bigl|\log|f_\lambda'(\psi_\lambda(\sigma^j\omega))|
       -\log|f_\lambda'(\psi_\lambda(\sigma^j\omega'))|\bigr|
 \le \sum_{j=0}^{n-1}C\theta^{\beta(n-j)}
 \le \frac{C\theta^\beta}{1-\theta^\beta},
\end{gather*}
so that for every $\lambda\in U_0$, every $n\in\bN$, and every $t\in\bR$, 
noting that every root $z$ of $f_\lambda^n(z)=z$ is written as 
$\psi_\lambda(\overline{k_0,k_1,k_2,\ldots,k_{n-1}})$ for some 
$(k_0,k_1,\ldots,k_{n-1})\in\{0,\ldots,d-1\}^n$,
we compute
\begin{multline*}
\int_{\bC}\bigl|(f_\lambda^n)'\bigr|^{-t}[f_\lambda^n-\Id=0]
 \le \exp\Bigl(\frac{|t|C\theta^\beta}{1-\theta^\beta}\Bigr)(\cL_{\lambda,t}^n1)(\psi_\lambda(\overline{0}))=\exp\Bigl(\frac{|t|C\theta^\beta}{1-\theta^\beta}\Bigr)
((\psi_\lambda^*\circ\cL_{\lambda,t}\circ(\psi_\lambda)_*)^n1)(\overline{0})\\
 \le \exp\Bigl(\frac{|t| C\theta^\beta}{1-\theta^\beta}\Bigr)
 \cdot\frac{1}{\inf_{\Sigma_d}h^\Sigma_{\lambda,t}}\Lambda_{\lambda,t}^n\cdot h^\Sigma_{\lambda,t}(\overline{0})
 \le \exp\Bigl(\frac{|t| C\theta^\beta}{1-\theta^\beta}\Bigr)\cdot
 \frac{\sup_{\Sigma_d}h^\Sigma_{\lambda,t}}{\inf_{\Sigma_d}h^\Sigma_{\lambda,t}}\Lambda_{\lambda,t}^n.
\end{multline*}
By the above joint continuity of $h^\Sigma_{\lambda,t}$ in $(\lambda,t)$
with values in $C^\alpha(\Sigma_d)$, the ratio
\begin{gather*}
 \frac{\sup_{\Sigma_d}h^\Sigma_{\lambda,t}}{\inf_{\Sigma_d}h^\Sigma_{\lambda,t}}
 =\frac{\sup_{J_{f_\lambda}}h_{\lambda,t}}{\inf_{J_{f_\lambda}}h_{\lambda,t}}
\end{gather*}
is continuous in $(\lambda,t)$ near $(\lambda_0,t_0)$ for any $t_0\in\bR$.
In particular, taking $t_0=0$, by a compactness argument,
we have the desired uniform order estimate \eqref{eq:sum}.

\section{An ergodic proof of Lemma \ref{th:upper}}\label{sec:ergodic}

We show Lemma \ref{th:upper} by a $d$-symbol version of 
the ergodic argument in \cite[Theorem 4.5 and \S 4.5]{FG21}, 
in a qualitative way (i.e., with no non-trivial order estimate).
We retain the setting and the notations in Appendix \ref{th:thermo}.
The uniform Bernoulli probability measure on $\Sigma_d$ is denoted by $m_d$.

\begin{lemma}\label{th:average}
Locally uniformly in $\lambda\in E_d$,
\begin{gather}\label{eq:average}
 \lim_{n\to\infty}\frac{1}{\nu(n)}\int_{\bC}
 \biggl|\frac{1}{n}\log\bigl|(f_\lambda^n)'\bigr|-L_d(\lambda)\biggr|
 [\Phi_{f,n}^*(\lambda,\cdot)=0]=0.
\end{gather}
\end{lemma}

\begin{proof}
Recall $\lambda_0\in E_d$, $U_0=U_0(\lambda_0,r_0)=\{\|\lambda-\lambda_0\|<r_0\}$,
$0<r_0\ll 1$, and $(\psi_\lambda)_{\lambda\in U_0}$ as in Appendix \ref{th:thermo}. 
For every $\lambda\in U_0$ and every $n\in\bN$, we compute
\begin{align*}
 &\frac{1}{d^n}\int_{\bC}\biggl|\frac{1}{n}\log\bigl|(f_\lambda^n)'\bigr|-L_d(\lambda)\biggr|[f_\lambda^n-\Id=0]\\
 =&\int_{\Sigma_d}
 \biggl|\frac{1}{n}\log\bigl|(f_\lambda^n)'
 (\psi_\lambda(\overline{\omega_0,\ldots,\omega_{n-1}}))\bigr|
 -L_d(\lambda)\biggr|\,m_d\\
 \le& \int_{\Sigma_d}
 \biggl|\frac{1}{n}\sum_{j=0}^{n-1}
 \log\bigl|f_\lambda'(\psi_\lambda(\sigma^j\omega))\bigr|
 -L_d(\lambda)\biggr|\,m_d\\
 &+\frac{1}{n}\cdot n\sup_{\omega\in\Sigma_d}
 \biggl|\frac{1}{n}\log\bigl|(f_\lambda^n)'
 (\psi_\lambda(\overline{\omega_0,\ldots,\omega_{n-1}}))\bigr|
 -\frac{1}{n}\sum_{j=0}^{n-1}\log\bigl|f_\lambda'(\psi_\lambda(\sigma^j\omega))\bigr|\biggr|\\
\le& \int_{\Sigma_d}
 \biggl|\frac{1}{n}\sum_{j=0}^{n-1}
 \log\bigl|f_\lambda'(\psi_\lambda(\sigma^j\omega))\bigr|
 -L_d(\lambda)\biggr|\,m_d\\
 &+
\frac{1}{n}\sum_{r=1}^n
 \sup_{\lambda\in U_0}
 \Bigl(\sup_{\omega,\omega'\in\Sigma_d:\, d_\theta(\omega,\omega')\le \theta^r}
 \bigl|\log|f_\lambda'(\psi_\lambda(\omega))|
 -\log|f_\lambda'(\psi_\lambda(\omega'))|\bigr|\Bigr),
\end{align*}
and by the uniqueness of $\mu_{f_\lambda}$,
we have $(\psi_\lambda)^*\mu_{f_\lambda}=m_d$ on $\Sigma_d$, so that
by Birkhoff's ergodic theorem and Fatou's lemma, we have
\begin{gather*}
 \int_{\Sigma_d}
 \biggl|\frac{1}{n}\sum_{j=0}^{n-1}
 \log\bigl|f_\lambda'(\psi_\lambda(\sigma^j\omega))\bigr|
 -L_d(\lambda)\biggr|\,m_d\to 0
 \quad\text{as }n\to\infty.
\end{gather*}
Moreover, as $0<r_0\searrow0$,
\begin{multline*}
 \sup_{\lambda,\lambda'\in U_0}\Biggl|\int_{\Sigma_d}
 \biggl|\frac{1}{n}\sum_{j=0}^{n-1}
 \log\bigl|f_\lambda'(\psi_\lambda(\sigma^j\omega))\bigr|
 -L_d(\lambda)\biggr|\,m_d
 -\int_{\Sigma_d}
 \biggl|\frac{1}{n}\sum_{j=0}^{n-1}
 \log\bigl|f_{\lambda'}'(\psi_{\lambda'}(\sigma^j\omega))\bigr|
 -L_d(\lambda')\biggr|\,m_d\Biggr|\\
 \le \sup_{\lambda,\lambda'\in U_0}\bigl\|\log|f_\lambda'\circ\psi_\lambda|
 -\log|f_{\lambda'}'\circ\psi_{\lambda'}|\bigr\|_{C^0(\Sigma_d)}
+\sup_{\lambda,\lambda'\in U_0}|L_d(\lambda)-L_d(\lambda')|\to 0.
\end{multline*}
Pick any $\epsilon>0$. Choosing $0<r_0\ll 1$ for each $\lambda_0\in E_d$ and,
for every compact subset $E$ in $E_d$, choosing a finite subcover of $E$
from $(U_0(\lambda_0,r_0))_{\lambda_0\in E}$, also by \eqref{eq:bound}, we have
\begin{gather*}
 \limsup_{n\to\infty}\sup_{\lambda\in E}\int_{\Sigma_d}
 \biggl|\frac{1}{n}\sum_{j=0}^{n-1}
 \log\bigl|f_\lambda'(\psi_\lambda(\sigma^j\omega))\bigr|
 -L_d(\lambda)\biggr|\,m_d<\epsilon,
\end{gather*}
and we are done by $[\Phi_{f,n}^*(\lambda,\cdot)=0]\le [f_\lambda^n-\Id=0]$
and $\nu(n)\sim d^n$.
\end{proof}

\begin{corollary}\label{th:escape}
For every compact subset $E$ in $E_d$ and
every $(w_n)_{n\in\bN}$ in $\bC$, if 
\begin{gather*}
 \limsup_{n\to\infty}\log|w_n|^{1/n}
 <\inf_{\lambda\in E}L_d(\lambda)
\end{gather*}
under the convention that $\log 0=-\infty$, 
then 
\begin{gather*}
 \limsup_{n\to\infty}\sup_{\lambda\in E}\frac{1}{\nu(n)}
 \log\biggl|\frac{p_{f,n}^*(\lambda,w_n)}{p_{f,n}^*(\lambda,0)}\biggr|
 \le 0.
\end{gather*}
\end{corollary}

\begin{proof}
Using $\log|1-z|\le \log 2+\max\{0,\log|z|\}$ and \eqref{eq:diff}, 
if $n\gg 1$, then for every $\lambda\in E$,
\begin{align*}
 \frac{1}{\nu(n)}
 \log\biggl|\frac{p_{f,n}^*(\lambda,w_n)}{p_{f,n}^*(\lambda,0)}\biggr|
 &\le \biggl(\frac{\log 2}{n}
 +\frac{1}{\nu(n)}\int_{\bC}
 \max\biggl\{0,\frac{\log|w_n|-\log|(f_\lambda^n)'|}{n}\biggr\}
 [\Phi_{f,n}^*(\lambda,\cdot)=0]\biggr)\\
 &\le \biggl(\frac{\log 2}{n}
 +\frac{1}{\nu(n)}\int_{\bC}
 \max\biggl\{0,L_d(\lambda)-\frac{1}{n}\log|(f_\lambda^n)'|\biggr\}
 [\Phi_{f,n}^*(\lambda,\cdot)=0]\biggr)\\
 &\le \biggl(\frac{\log 2}{n}
 +\frac{1}{\nu(n)}\int_{\bC}
 \biggl|\frac{1}{n}\log\bigl|(f_\lambda^n)'\bigr|-L_d(\lambda)\biggr|
 [\Phi_{f,n}^*(\lambda,\cdot)=0]\biggr),
\end{align*}
and we are done by Lemma \ref{th:average}.
\end{proof}

\begin{proof}[Proof of Lemma \ref{th:upper} using Corollary \ref{th:escape}]
On every compact
subset $E$ in $E_d\cap\Omega_d(\ell)$, by the definition of $\Omega_d(\ell)$ and $\log|a_n|/n\to\ell$,
\begin{gather*}
 \limsup_{n\to\infty}\frac{\log|w_n|}{n}
 =\log d+\ell<\inf_{\lambda\in E}L_d(\lambda),
\end{gather*}
and we are done by Corollary \ref{th:escape}.
\end{proof}

\begin{acknowledgement}
 This research was partially supported by JSPS Grant-in-Aid for Scientific Research (C), 23K03129.
\end{acknowledgement}


\begin{thebibliography}{10}

\bibitem{BassanelliBerteloot09}
Giovanni Bassanelli and Fran{\c{c}}ois Berteloot.
\newblock Lyapunov exponents, bifurcation currents and laminations in
  bifurcation loci.
\newblock {\em Mathematische Annalen}, 345(1):1--23, 2009.

\bibitem{BassanelliBerteloot11}
Giovanni Bassanelli and Fran{\c{c}}ois Berteloot.
\newblock Distribution of polynomials with cycles of a given multiplier.
\newblock {\em Nagoya Math. J.}, 201:23--43, 2011.

\bibitem{BM01}
Fran{\c{c}}ois Berteloot and Volker Mayer.
\newblock {\em Rudiments de dynamique holomorphe}, volume~7 of {\em Cours
  Sp\'ecialis\'es [Specialized Courses]}.
\newblock Soci\'et\'e Math\'ematique de France, Paris, 2001.

\bibitem{BrannerHubbard88}
Bodil Branner and John~H. Hubbard.
\newblock The iteration of cubic polynomials. {Part I}: The global topology of
  parameter space.
\newblock {\em Acta Math.}, 160:143--206, 1988.

\bibitem{BrannerHubbard92}
Bodil Branner and John~H. Hubbard.
\newblock The iteration of cubic polynomials. {Part II}: Patterns and parapatterns.
\newblock {\em Acta Math.}, 169(3--4):229--325, 1992.

\bibitem{Brolin}
Hans Brolin.
\newblock Invariant sets under iteration of rational functions.
\newblock {\em Ark. Mat.}, 6:103--144, 1965.

\bibitem{BG14}
Xavier Buff and Thomas Gauthier.
\newblock Quadratic polynomials, multipliers and equidistribution.
\newblock {\em Proc. Amer. Math. Soc.}, 143(7):3011--3017, 2015.

\bibitem{DH84}
Adrien Douady and John~H. Hubbard.
\newblock {\em \'Etude dynamique des polyn\^omes complexes. Premi\`ere partie}.
\newblock Publications Math\'ematiques d'Orsay, 84-02, Universit\'e de Paris-Sud, 1984.

\bibitem{DF08}
Romain Dujardin and Charles Favre.
\newblock Distribution of rational maps with a preperiodic critical point.
\newblock {\em Amer. J. Math.}, 130:979--1032, 2008.

\bibitem{EL97}
AE~Er{\"e}menko and GM~Levin.
\newblock Estimating the characteristic exponents of polynomials.
\newblock {\em Journal of Mathematical Sciences}, 85(5):2164--2171, 1997.

\bibitem{FG21}
Tanya Firsova and Igors Gorbovickis.
\newblock Equidistribution of critical points of the multipliers in the
  quadratic family.
\newblock {\em Adv. Math.}, 380:107591, 2021.

\bibitem{Gauthier13}
Thomas Gauthier.
\newblock Equidistribution towards the bifurcation current {I}: multipliers and
  degree {$d$} polynomials.
\newblock {\em Math. Ann.}, 366(1-2):1--30, 2016.

\bibitem{GauthierOkuyamaVigny19}
Thomas Gauthier, Y\^usuke Okuyama, and Gabriel Vigny.
\newblock Hyperbolic components of rational maps: Quantitative equidistribution
  and counting.
\newblock {\em Comment. Math. Helv.}, 94:347--398, 2019.

\bibitem{GauthierVigny17}
Thomas Gauthier and Gabriel Vigny.
\newblock Distribution of postcritically finite polynomials {II}: Speed of
  convergence.
\newblock {\em J. Mod. Dyn.}, 11:57--98, 2017.

\bibitem{Goldberg94}
Lisa~R. Goldberg.
\newblock On the multiplier of a repelling fixed point.
\newblock {\em Invent. Math.}, 118:85--108, 1994.

\bibitem{GS00}
Jacek Graczyk and Grzegorz {\'S}wi{\c{a}}tek.
\newblock Harmonic measure and expansion on the boundary of the connectedness
  locus.
\newblock {\em Inventiones mathematicae}, 142(3):605--629, 2000.

\bibitem{GS17}
Jacek Graczyk and Grzegorz {\'S}wi{\k{a}}tek.
\newblock Fine structure of connectedness loci.
\newblock {\em Math. Ann.}, 369(1-2):49--108, 2017.

\bibitem{Hormander83}
Lars H{\"o}rmander.
\newblock {\em The analysis of linear partial differential operators. {I}},
  volume 256 of {\em Grundlehren der Mathematischen Wissenschaften [Fundamental
  Principles of Mathematical Sciences]}.
\newblock Springer-Verlag, Berlin, 1983.
\newblock Distribution theory and Fourier analysis.

\bibitem{ImayoshiTaniguchi92}
Yoichi Imayoshi and Masahiko Taniguchi.
\newblock {\em An Introduction to Teichm{\"u}ller Spaces}.
\newblock Springer-Verlag, Tokyo, 1992.

\bibitem{Kiwi97}
Jan Kiwi.
\newblock Rational rays and critical portraits of complex polynomials.
\newblock Stony Brook IMS Preprint 1997/15, arXiv:math/9710212v1, 1997.

\bibitem{Kiwi04}
Jan Kiwi.
\newblock Real laminations and the topological dynamics of complex polynomials.
\newblock {\em Adv. Math.}, 184(2):207--267, 2004.

\bibitem{Levin90}
Gennadi~M Levin.
\newblock Theory of iterations of polynomial families in the complex plane.
\newblock {\em Journal of Mathematical Sciences}, 52(6):3512--3522, 1990.

\bibitem{MSS}
R.~Ma{\~n}\'e, P.~Sad, and D.~Sullivan.
\newblock On the dynamics of rational maps.
\newblock {\em Ann. Sc. E.N.S.,4$\grave{e}$me s$\acute{e}$rie}, 16:193--217,
  1983.

\bibitem{Manning84}
Anthony Manning.
\newblock The dimension of the maximal measure for a polynomial map.
\newblock {\em Ann. of Math. (2)}, 119(2):425--430, 1984.

\bibitem{McMullen94}
Curtis~T. McMullen.
\newblock {\em Complex Dynamics and Renormalization}, volume~135 of
  {\em Annals of Mathematics Studies}.
\newblock Princeton University Press, Princeton, 1994.

\bibitem{MortonSilverman95}
Patrick Morton and Joseph~H. Silverman.
\newblock Periodic points, multiplicities, and dynamical units.
\newblock {\em J. Reine Angew. Math.}, 461:81--122, 1995.

\bibitem{MortonVivaldi95}
Patrick Morton and Franco Vivaldi.
\newblock Bifurcations and discriminants for polynomial maps.
\newblock {\em Nonlinearity}, 8(4):571--584, 1995.

\bibitem{Okuyama18}
Y\^usuke Okuyama.
\newblock Nevanlinna theory and value distribution in the unicritical
  polynomials family.
\newblock {\em Riv. Mat. Univ. Parma (N.S.)}, 9(1):1--19, 2018.

\bibitem{Przytycki85}
Feliks Przytycki.
\newblock Hausdorff dimension of harmonic measure on the boundary of an
  attractive basin for a holomorphic map.
\newblock {\em Invent. Math.}, 80(1):161--179, 1985.

\bibitem{PU10}
Feliks Przytycki and Mariusz Urba{\'n}ski.
\newblock {\em Conformal Fractals: Ergodic Theory Methods}.
\newblock Cambridge University Press, Cambridge, 2010.

\bibitem{Schleicher99}
Dierk Schleicher.
\newblock On fibers and local connectivity of compact sets in {$\bC$}.
\newblock Stony Brook IMS Preprint 1998/12, 1998.
\newblock arXiv:math/9902154v1.

\bibitem{Shishikura}
Mitsuhiro Shishikura.
\newblock On the quasiconformal surgery of rational functions.
\newblock {\em Ann. Sci. \'{E}cole Norm. Sup. (4)}, 20(1):1--29, 1987.

\bibitem{Sullivan85}
Dennis Sullivan.
\newblock Quasiconformal homeomorphisms and dynamics. {I}. Solution of the
  {Fatou--Julia} problem on wandering domains.
\newblock {\em Ann. of Math. (2)}, 122:401--418, 1985.

\bibitem{Zakeri06}
Saeed Zakeri.
\newblock On biaccessible points of the Mandelbrot set.
\newblock {\em Proceedings of the American Mathematical Society},
  134(8):2239--2250, 2006.

\bibitem{Zins00}
Michel Zinsmeister.
\newblock {\em Thermodynamic formalism and holomorphic dynamical systems},
  volume~4.
\newblock American Mathematical Soc., 2000.

\end{thebibliography}
\end{document}